\PassOptionsToPackage{reqno}{amsmath}
\documentclass[hidelinks]{siamart251216}
\usepackage[T1]{fontenc}
\usepackage[utf8]{inputenc}
\usepackage{microtype}
\usepackage{amssymb,mathtools,mathrsfs,bm}
\usepackage[ntheorem]{empheq}
\usepackage{etoolbox}
\makeatletter
\MHInternalSyntaxOn
\patchcmd{\EQ_post_process:}
{\Hy@linkcounter=\EQ_hyperref_link_fint}
{}{}%
{\PackageError{ALE-SAV-ContactLine}{Cannot patch empheq hyperlink counter}{}}
\MHInternalSyntaxOff
\makeatother
\usepackage{booktabs}
\usepackage{graphicx}
\makeatletter
\let\siam@makecaption\@makecaption
\makeatother
\usepackage{subcaption}
\makeatletter
\AtBeginDocument{\let\@makecaption\siam@makecaption}
\makeatother
\usepackage{xcolor}
\usepackage{placeins}
\usepackage[numbers,sort&compress]{natbib}
\makeatletter
\pretocmd{\NAT@citenum}{\color{black}}{}%
{\PackageError{ALE-SAV-ContactLine}{Cannot set citation colour}{}}
\makeatother
\hypersetup{pdftitle={A linearly implicit high-order structure-preserving ALE-SAV finite element method for moving contact line flows},pdfauthor={Shu Ma}}
\allowdisplaybreaks
\AtBeginEnvironment{thebibliography}{\color{black}\interlinepenalty=10000\relax}
\newsiamthm{assumption}{Assumption}
\newsiamremark{remark}{Remark}
\newcommand{\R}{\mathbb R}
\newcommand{\Om}{\Omega}
\newcommand{\Ga}{\Gamma}
\newcommand{\D}{\bm D}
\newcommand{\dd}{\,\mathrm d}

\newcommand{\Div}{\nabla\!\cdot}

\newcommand{\E}{\mathcal E}   
\newcommand{\V}{\mathcal V}
\newcommand{\MDR}{\mathcal M}

\newcommand{\norm}[1]{\left\lVert#1\right\rVert}

\makeatletter
\newcommand{\algorithmstep}{\@startsection{paragraph}{4}{0pt}%
{-1.3ex plus -.5ex minus -.2ex}{-.5em plus -.1em}%
{\normalfont\normalsize\bfseries}}
\makeatother
\title{A high-order structure-preserving ALE--SAV finite element method for moving contact line flows}
\author{Shu Ma\thanks{Department of Mathematics, Hong Kong Baptist University, Kowloon Tong, Hong Kong (\email{shu\_ma@hkbu.edu.hk}).}}
\headers{High-order structure-preserving ALE--SAV method}{Shu Ma}
\begin{document}
\maketitle

\begin{abstract}
We develop a structure-preserving arbitrary Lagrangian--Eulerian (ALE)
finite element method for incompressible
Navier--Stokes flow with gravity and moving contact lines, based on a
scalar auxiliary variable (SAV) reformulation. The reformulation
expresses the geometric energy quadratically in the auxiliary variable
and permits linear treatment of the capillary, wetting and gravitational
forces. The method combines high-order isoparametric elements in space
with a linearly implicit, second-order Crank--Nicolson discretization in
time. Second-order extrapolation supplies the predicted geometry and
known coefficients of the linear fluid--SAV system. Mesh motion changes
the kinetic-energy inner product and requires a compatible discretization
of the time-derivative terms in the ALE formulation.
We discretize these terms using an
$L^2$-isometric velocity transfer to retain the kinetic-energy balance
after the geometry update. We also construct a boundary update that
preserves area or volume exactly using only linear solves. The resulting
scheme satisfies a discrete SAV energy-dissipation
law and preserves area or volume exactly. Numerical experiments in two
and three dimensions demonstrate high-order spatial and second-order
temporal convergence and confirm SAV energy dissipation and area or
volume conservation.
\end{abstract}


\section{Introduction}\label{sec:physical-model}\label{sec:physical-energy}
\begin{figure}[!t]
\centering
\includegraphics[width=.55\linewidth]{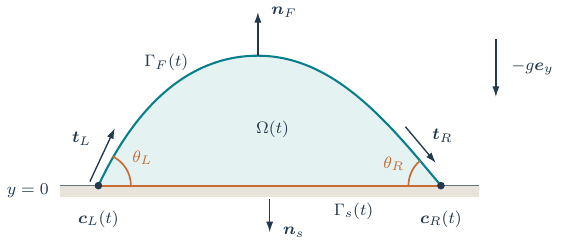} \hspace{-10pt}  
\resizebox{.55\linewidth}{!}{%
\includegraphics[trim=0 14bp 0 0,clip]{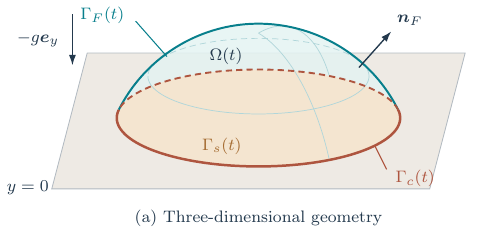}%
\includegraphics[trim=0 14bp 0 0,clip]{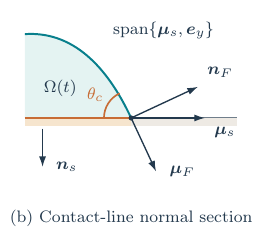}%
}
\caption{Droplet geometry in two dimensions (top) and three dimensions (bottom).}
\label{fig:model-geometry}
\label{fig:model-geometry-3d}
% \vspace{-10pt}
\end{figure}

Let $\Om(t)\subset\{y>0\}$ be a bounded 
domain in $d\in\{2,3\}$, where $y=x_d$.
Its free surface $\Ga_F(t)$ and wetted wall
$\Ga_s(t)\subset\{y=0\}$ meet at
$\Ga_c(t)=\partial\Ga_F(t)=\partial\Ga_s(t)$. 
The free surface $\Ga_F(t)$ has outward unit normal $\bm n_F$
and normal speed $V_F$. We use the curvature convention
$\kappa\bm n_F=-\Delta_{\Ga_F}\operatorname{id}$.
The wall normal is $\bm n_s=-\bm e_y$, where $\bm e_y$ is the
vertical unit vector, and $P_s=I-\bm n_s\otimes\bm n_s$ is the
tangential projection.
The geometry is illustrated in Figure~\ref{fig:model-geometry}.
The flow is governed by the incompressible Navier--Stokes equations
with gravity and moving-contact-line boundary conditions
\cite{HLLT,GHL}:
\begin{subequations}\label{eq:physical-system}
\begin{empheq}[left=\empheqlbrace]{alignat=3}
&\text{\footnotesize\upshape (fluid)}&\;&\partial_t\bm u+(\bm u\cdot\nabla)\bm u-\Div\bm\sigma
=-g\bm e_y&\;&\text{in }\Om(t),\label{eq:NS}\\
&&\;&\Div \bm u=0&\;&\text{in }\Om(t).\label{eq:incompressibility}\\[.8ex]
&\text{\footnotesize\upshape (free surface)}&\;&\bm\sigma \bm n_F=-\gamma\kappa \bm n_F&\;&\text{on }\Ga_F(t),\label{eq:free-traction}\\
&&\;&V_F=\bm u\cdot\bm n_F&\;&\text{on }\Ga_F(t).\label{eq:free-boundary-continuous}\\[.8ex]
&\text{\footnotesize\upshape (wall)}&\;&\bm u\cdot \bm n_s=0&\;&\text{on }\Ga_s(t),\label{eq:wall-conditions}\\
&&\;&P_s(\bm\sigma \bm n_s+\beta_s\bm u)=0&\;&\text{on }\Ga_s(t),\label{eq:wall-slip}\\[.8ex]
&\text{\footnotesize\upshape (contact line)}&\;&\beta_cU_c=\gamma(\cos\theta_Y-\cos\theta_c)&\;&\text{on }\Ga_c(t),\label{eq:younglaw-general}
\end{empheq}
\end{subequations}
where  the  stress is
defined by
\begin{equation}\label{eq:physical-stress-definition}
\D\bm u=\tfrac12(\nabla\bm u+\nabla\bm u^T),\qquad
\bm\sigma=2\mu\D\bm u-p_{\rm phys}I.
\end{equation}
Along $\Ga_c(t)$, let $\bm\mu_s$ and $\bm\mu_F$ be the outward
unit conormals of the wetted wall within the substrate and of
the free surface, respectively. The liquid contact angle
$\theta_c\in(0,\pi)$ satisfies
$\cos\theta_c=\bm\mu_F\cdot\bm\mu_s$.
The kinematic condition \eqref{eq:free-boundary-continuous}
and the wall impermeability condition \eqref{eq:wall-conditions}
give $U_c=\bm u\cdot\bm\mu_s$, the contact speed along $\bm\mu_s$.

The viscous and surface-tension coefficients are
$\mu=\mathrm{Re}^{-1}>0$ and $\gamma=\mathrm{We}^{-1}>0$
(inverse Reynolds and Weber numbers).
The coefficients $\beta_s,\beta_c\geq0$ describe wall slip
and contact-line friction, respectively.
Gravity is $-g\bm e_y$ with $g\geq0$.
Young's equilibrium angle is $\theta_Y\in(0,\pi)$.
For $\beta_c=0$, \eqref{eq:younglaw-general} gives
$\theta_c=\theta_Y$.

 Let $\bm X(t)$ parametrize $\Ga_F(t)$.
The initial data are $\bm X(\cdot,0)=\bm X_0$ and
$\bm u(\cdot,0)=\bm u_0$. The geometric energy $\E_{\rm G}$ and total physical energy $\E$ are
\begin{align}
\E_{\rm G}(\bm X)&=\gamma|\Ga_F(t)|-\gamma\cos\theta_Y\,|\Ga_s(t)|
+ g\int_{\Om(t)}y\dd x,\label{eq:geometric-energy}\\
\E(t)&=\frac12\norm{\bm u}_{L^2(\Om(t))}^2
+\E_{\rm G}(\bm X(t)).\label{eq:physical-total-energy}
\end{align}
For sufficiently regular solutions, testing \eqref{eq:NS} by $\bm u$,
applying Reynolds' transport theorem and using \eqref{eq:free-traction},
\eqref{eq:wall-slip} and \eqref{eq:younglaw-general} gives \cite{GHL}
\begin{equation}\label{eq:continuous-energy}
\frac{\dd}{\dd t}\E(t)
+2\mu\norm{\D\bm u}_{L^2(\Om(t))}^2
+\beta_s\norm{P_s\bm u}_{L^2(\Ga_s(t))}^2
+\beta_c\int_{\Ga_c(t)}|\bm u\cdot\bm\mu_s|^2\dd\ell=0.
\end{equation}
Here and throughout, $\dd\ell$ denotes counting measure
(endpoint summation) for $d=2$ and arclength measure for $d=3$.
Writing $\V(\bm X)=|\Om(t)|$, we also obtain from
\eqref{eq:free-boundary-continuous}, \eqref{eq:incompressibility} and
\eqref{eq:wall-conditions} 
\begin{equation}\label{eq:continuous-volume-conservation}
\begin{aligned}
\frac{\dd}{\dd t}\V(\bm X(t))
&=\int_{\Ga_F(t)}V_F\dd s
=\int_{\Ga_F(t)}\bm u\cdot\bm n_F\dd s =\int_{\Om(t)}\Div\bm u\dd x=0.
\end{aligned}
\end{equation}
The divergence theorem gives $|\Ga_s(t)|\leq|\Ga_F(t)|$,
so $g\geq0$ and $y\geq0$ imply
\begin{equation}\label{eq:geometric-energy-lower-bound}
\E_{\rm G}(\bm X)\geq\gamma(1-\cos\theta_Y)|\Ga_s(t)|\geq0.
\end{equation}

For finite element approximation on evolving fluid domains, an
arbitrary Lagrangian--Eulerian (ALE) formulation allows the interior
mesh velocity to be chosen independently of the fluid velocity
\cite{FormaggiaNobileALE2004,DuanLiYang,GarckeNurnbergZhao,
GarckeNurnbergZhaoAxisymmetricALE,LiMaQiuALE}.
Accuracy and structure preservation for evolving interfaces have also been studied using parametric finite element methods. 
High-order time discretizations of the Barrett--Garcke--N\"urnberg
(BGN) formulation have been developed in
\cite{JiangSuZhangBGN2,JiangSuZhangBDF,MaRaoBGN}, and high-order
isoparametric approximations have been studied in
\cite{GarckeNurnbergPraetoriusZhang}.
Bao and Zhao \cite{BaoZhaoSurfaceDiffusion} constructed a parametric
finite element method for surface diffusion that preserves volume
exactly and dissipates geometric energy.
Nonlinear Lagrange multiplier formulations also preserve volume and
satisfy a geometric energy-dissipation law \cite{GarckeJiangSuZhangLM}. 
% Discrete structure preservation remains a key issue in the numerical
% approximation of evolving interfaces. 
For the coupled fluid problem with moving contact lines,
energy-stable sharp-interface discretizations incorporate wall slip and
dynamic contact angles \cite{ZhaoRen,IS}.
The first-order ALE methods of Hu et al.\ \cite{HLLT} and Garcke,
Hu and Lei \cite{GHL} satisfy a discrete dissipation law for the physical
energy.
Their exact-volume variants require nonlinear solves.

Combining high-order approximation with linear solves while retaining a discrete energy law and exact volume conservation remains challenging for moving-contact-line flows. Mesh motion changes the kinetic-energy inner product.
In the skew-symmetric ALE formulation
\eqref{eq:ale-skew-equivalence}, the time-derivative terms are
$\partial_t^\circ\bm u+\tfrac12(\Div\bm w)\bm u$,
where $\bm w$ is the mesh velocity and $\partial_t^\circ$
is ALE derivative defined in \eqref{eq:ale-continuous}.
These terms require a compatible time discretization that retains
the kinetic-energy balance on the final mesh. In addition, discrete incompressibility and wall impermeability on the predicted mesh provide a zero-flux condition, but do not by themselves ensure volume conservation on the final mesh. The kinetic-energy discretization and the volume-preserving boundary update must therefore be compatible within a high-order scheme requiring only linear solves.

% To treat the capillary, wetting and gravitational forces linearly in the discretization, we use the scalar auxiliary variable (SAV)
% approach, which introduces a quadratic modified energy
% \cite{ShenXuYang,ShenXuYangReview}. The SAV approach has been applied to
% high-order time stepping \cite{AkrivisLiLi,FengLiMaSAV}, incompressible
% Navier--Stokes and phase-field two-phase flows on fixed domains
% \cite{LiShenLiuSAVNS,LiShenMSAV}, problems with global constraints enforced through
% additional nonlinear solves \cite{ChengShenConstraints}, and ALE formulations for closed-interface flow
% \cite{MaSAV}.

% We construct an ALE--SAV finite element method for
% \eqref{eq:physical-system} using high-order isoparametric elements in
% space and a linearly implicit, second-order Crank--Nicolson
% discretization in time.
% The auxiliary variable accounts for the capillary, wetting and
% gravitational energies, while wall slip and contact-line friction
% enter the dissipation.
% Second-order extrapolation supplies the predicted geometry and known
% coefficients of the linear fluid--SAV system. The following two constructions yield a discrete SAV energy-dissipation law and exact area or volume conservation.

We develop a high-order, linearly implicit ALE-SAV finite element method for \eqref{eq:physical-system} that satisfies a discrete modified-energy dissipation law and preserves volume exactly. To treat the capillary, wetting and gravitational forces linearly in the discretization, we use the scalar auxiliary variable (SAV) approach, which introduces a quadratic modified energy \cite{ShenXuYang,ShenXuYangReview}. The SAV approach has been applied to high-order time stepping \cite{AkrivisLiLi,FengLiMaSAV}, incompressible Navier-Stokes and phase-field two-phase flows on fixed domains \cite{LiShenLiuSAVNS,LiShenMSAV} and  closed-interface flow \cite{MaSAV}. The ALE-SAV method combines high-order isoparametric elements in space with a linearly implicit, second-order Crank-Nicolson discretization in time. Wall slip and contact-line friction enter the dissipation. Second-order extrapolation supplies the predicted geometry and known coefficients of the linear fluid-SAV system. The following two constructions address the kinetic-energy balance on moving meshes and the volume constraint.

\smallskip

\setlength{\leftmargini}{1.9em}
\setlength{\labelwidth}{1.25em}
\setlength{\labelsep}{0.35em}
\renewcommand{\labelenumi}{(\arabic{enumi})}
\begin{enumerate}

\item \textit{Kinetic-energy balance on moving meshes.}
We construct an \(L^2\)-isometric velocity transfer
between evolving finite element spaces and use it, together with
its associated skew-symmetric correction, to discretize these
time-derivative terms. The endpoint velocity update retains the kinetic-energy balance on the final mesh. Together with the linear SAV coupling, this balance yields a discrete modified-energy dissipation law.
\smallskip

\item \textit{Exact area or volume conservation through linear solves.}
We construct a boundary map from midpoint substeps that preserve quadratic invariants \cite{SMB}, with frozen skew-symmetric operators. Area is quadratic in two dimensions. In three dimensions, a symmetric composition updates two coordinate groups at each substep, so volume is at most quadratic in the active variables. Each substep is therefore linear and preserves area or volume exactly.
\end{enumerate}
\smallskip

The rest of the paper is organized as follows. 
Section~\ref{sec:model} introduces the ALE--SAV
reformulation. Section~\ref{sec:space} defines the evolving finite
element spaces and structure-preserving operators.
Section~\ref{sec:scheme} presents the fully discrete scheme. Section~\ref{sec:proofs} establishes the discrete
energy and volume laws.
Section~\ref{sec:diagnostic} presents numerical experiments in two and
three dimensions.

% \new page
\section{The ALE--SAV reformulation}\label{sec:model}
In this section, we reformulate \eqref{eq:physical-system} in an ALE--SAV framework
and derive its weak formulation and the corresponding volume
conservation and energy dissipation laws.

\subsection{ALE flow map}\label{sec:extension-continuous}

Let $\Phi(\cdot,t):\Om(0)\to\Om(t)$ be the ALE flow map.
On any interval where $\Phi$ is a smooth diffeomorphism, its
ALE velocity $\bm w$ satisfies
\begin{equation}\label{eq:ale-map-evolution}
\partial_t\Phi(\bm x_0,t)=\bm w(\Phi(\bm x_0,t),t),\qquad
\Phi(\bm x_0,0)=\bm x_0,\quad \bm x_0\in\Om(0).
\end{equation}
The free-surface parametrization
$\bm X(\cdot,t):\widehat\Ga_F\to\Ga_F(t)$ uses a fixed reference
surface $\widehat\Ga_F$, with oriented reference patches for $d=3$.
For $d=2$, take $\widehat\Ga_F=[0,1]$.
The boundary trace of $\Phi$ maps $\Ga_F(0)$ and $\Ga_s(0)$ onto
$\Ga_F(t)$ and $\Ga_s(t)$, respectively, and satisfies
$\bm X(\cdot,t)=\Phi(\cdot,t)\circ\bm X(\cdot,0)$.

The corresponding ALE derivative is
\begin{equation}\label{eq:ale-continuous}
\partial_t^\circ\bm u
=\bigl[\partial_t(\bm u\circ\Phi)\bigr]\circ\Phi^{-1}
=\partial_t\bm u+(\bm w\cdot\nabla)\bm u.
\end{equation}
For a wall-impermeable field $\bm v$, let
$\bm w=\MDR_{\bm X}[\bm v]$ be the generalized harmonic extension
of its free-surface trace associated with the symmetric-gradient form. This ALE harmonic extension satisfies
\begin{subequations}\label{eq:mdr-strong}\label{eq:mesh-extension-traces}
\begin{empheq}[left=\empheqlbrace]{alignat=2}
-\Div(\D\bm w)&=0&\qquad&\text{in }\Om(t),\label{eq:mdr-bulk-strong}\\
\bm w&=\bm v&\qquad&\text{on }\Ga_F(t),\label{eq:mdr-free-trace-strong}\\
\bm w\cdot\bm n_s&=0&\qquad&\text{on }\Ga_s(t),\label{eq:mdr-wall-normal-strong}\\
P_s(\D\bm w\,\bm n_s)&=0&\qquad&\text{on }\Ga_s(t).\label{eq:mdr-wall-tangential-strong}
\end{empheq}
\end{subequations}

We choose the material parametrization
$\partial_t\bm X=\bm u\circ\bm X$ and set
$\bm w=\MDR_{\bm X}[\bm u]$, so
\eqref{eq:mdr-free-trace-strong} gives $\bm w=\bm u$ on $\Ga_F(t)$.
No incompressibility constraint is imposed on $\bm w$.

\subsection{SAV formulation}\label{sec:reformulation}\label{sec:geometric-sav}

To express the geometric energy quadratically, choose $C_0>0$ and
introduce the scalar auxiliary variable \cite{ShenXuYang} by
\begin{equation}\label{eq:sav-continuous}
r=\sqrt{\E_{\rm G}(\bm X)+C_0}.
\end{equation}
The lower bound \eqref{eq:geometric-energy-lower-bound} ensures a positive
square root, and this transformation gives $\E_{\rm G}=r^2-C_0$.
For a smooth deformation field $\bm\xi$ with
$\bm\xi\cdot\bm n_s=0$ on $\Ga_s(t)$, let
$D\E_{\rm G}(\bm X)[\bm\xi]$ denote the derivative of
\eqref{eq:geometric-energy} along
$(\mathrm{Id}+\varepsilon\bm\xi)\circ\bm X$ at $\varepsilon=0$.
Then
\begin{equation}\label{eq:DF}
\begin{aligned}
D\E_{\rm G}(\bm X)[\bm{\xi}]
&=\gamma\int_{\Ga_F(t)}
\nabla_{\Ga_F}\!\cdot\bm{\xi}\dd s
-\gamma\cos\theta_Y\int_{\Ga_c(t)}
\bm{\xi}\cdot\bm\mu_s\dd\ell +g\int_{\Ga_F(t)}y\,\bm{\xi}\cdot\bm n_F\dd s\\
&=\int_{\Ga_F(t)}(\gamma\kappa+ gy)\bm n_F\cdot\bm{\xi}\dd s+\gamma\int_{\Ga_c(t)}(\cos\theta_c-\cos\theta_Y)
\bm{\xi}\cdot\bm\mu_s\dd\ell,
\end{aligned}
\end{equation}
where the second form follows from surface integration by parts and
$\bm\mu_F\cdot\bm\xi=\cos\theta_c\,\bm\mu_s\cdot\bm\xi$
on the contact set.

Introduce the shifted pressure $\pi=p_{\rm phys}+gy$ and stress
$\bm\sigma_\pi=2\mu\D\bm u-\pi I=\bm\sigma-gyI$. Then
$\Div\bm\sigma_\pi=\Div\bm\sigma-g\bm e_y$ in the bulk, and
$gy$ enters the free-surface traction. Define the normalized first
variation by
\begin{equation}\label{eq:sav-normalized-force}
\mathcal G(\bm X)[\bm v]
=\frac{D\E_{\rm G}(\bm X)[\bm v]}
{\sqrt{\E_{\rm G}(\bm X)+C_0}}.
\end{equation}

We treat $r$ as an independent unknown whose evolution
equation is derived from \eqref{eq:sav-continuous} by the chain rule.
Using \eqref{eq:ale-map-evolution} and $\bm w=\MDR_{\bm X}[\bm u]$ in \eqref{eq:mdr-strong}, the coupled SAV system reads
\begin{subequations}\label{eq:ale-sav-strong}
\begin{empheq}[left=\empheqlbrace]{alignat=3}
&\text{\footnotesize\upshape (fluid)}&\;&\partial_t^\circ\bm u+((\bm u-\bm w)\cdot\nabla)\bm u
=\Div\bm\sigma_\pi&\;&\text{in }\Om(t),\label{eq:ale-sav-bulk-strong}\\
&&\;&\Div \bm u=0&\;&\text{in }\Om(t).\label{eq:ale-sav-div-strong}\\[.8ex]
&\text{\footnotesize\upshape (SAV)}&\;&\dot r=\tfrac12\mathcal G(\bm X)[\bm u]&\;&t>0.\label{eq:ale-sav-scalar-strong}\\[.8ex]
&\text{\footnotesize\upshape (free surface)}&\;&\bm\sigma_\pi \bm n_F=-\frac{r}{\sqrt{\E_{\rm G}(\bm X)+C_0}}(\gamma\kappa+ gy)\bm n_F&\;&\text{on }\Ga_F(t),\label{eq:ale-sav-traction-strong}\\
&&\;&\partial_t\bm X=\bm u\circ \bm X&\;&\text{on }\widehat\Ga_F.\label{eq:ale-sav-material-strong}\\[.8ex]
&\text{\footnotesize\upshape (wall)}&\;&\bm u\cdot \bm n_s=0&\;&\text{on }\Ga_s(t),\label{eq:ale-sav-wall-strong}\\
&&\;&P_s(\bm\sigma_\pi \bm n_s+\beta_s\bm u)=0&\;&\text{on }\Ga_s(t).\label{eq:ale-sav-slip-strong}\\[.8ex]
&\text{\footnotesize\upshape (contact line)}&\;&\beta_cU_c=\frac{r}{\sqrt{\E_{\rm G}(\bm X)+C_0}}\gamma(\cos\theta_Y-\cos\theta_c)&\;&\text{on }\Ga_c(t).\label{eq:ale-sav-young-general}
\end{empheq}
\end{subequations}
We supplement the physical initial data with
\begin{equation}\label{eq:sav-initial-data}
r(0)=\sqrt{\E_{\rm G}(\bm X_0)+C_0}.
\end{equation}
For sufficiently regular solutions,
\eqref{eq:ale-sav-scalar-strong} and
\eqref{eq:ale-sav-material-strong} imply
\begin{equation}\label{eq:sav-continuous-compatibility}
\frac{\dd}{\dd t}\left(r-\sqrt{\E_{\rm G}(\bm X)+C_0}\right)=0.
\end{equation}
Together with \eqref{eq:sav-initial-data},
\eqref{eq:sav-continuous-compatibility} recovers
\eqref{eq:sav-continuous}, so
$r\mathcal G(\bm X)[\bm v]=D\E_{\rm G}(\bm X)[\bm v]$.
With $p_{\rm phys}=\pi-gy$, the system
\eqref{eq:ale-sav-strong} therefore recovers \eqref{eq:physical-system}
with the chosen material parametrization.

\subsection{Weak formulation and structural properties}\label{sec:sav-energy-structure} 
We retain the ALE map \eqref{eq:ale-map-evolution} and the material surface motion
\eqref{eq:ale-sav-material-strong}, with the initial data above.
The weak form of the ALE harmonic extension \eqref{eq:mdr-strong} is
\begin{equation}\label{eq:mdr-continuous-weak}
(\D\bm w,\D\bm z)_{\Om(t)}=0,
\end{equation}
for all $\bm z\in H^1(\Om(t))^d$ with
$\bm z=0$ on $\Ga_F(t)$ and
$\bm z\cdot\bm n_s=0$ on $\Ga_s(t)$.

For a sufficiently regular advecting field $\bm a$, define the skew-symmetric
convection form by
\begin{equation}\label{eq:skew-convection-continuous}
\begin{aligned}
b(t;\bm a,\bm u,\bm v)
&:=\frac12\Bigl([\bm a\cdot\nabla]\bm u,\bm v\Bigr)_{\Om(t)} -\frac12\Bigl([\bm a\cdot\nabla]\bm v,\bm u\Bigr)_{\Om(t)}.
\end{aligned}
\end{equation}
It is skew-symmetric in its last two arguments, so
$b(t;\bm a,\bm u,\bm u)=0$.
Let $\bm n$ denote the outward unit normal of $\partial\Om(t)$.
If $\Div\bm u=0$ and
$(\bm u-\bm w)\cdot\bm n=0$ on $\partial\Om(t)$,
then \eqref{eq:ale-continuous} and integration by parts give
\begin{equation}\label{eq:ale-skew-equivalence}
\begin{aligned}
\Bigl(\partial_t\bm u+(\bm u\cdot\nabla)\bm u,\bm v\Bigr)_{\Om(t)}
&=(\partial_t^\circ\bm u,\bm v)_{\Om(t)}
+\frac12\bigl((\Div\bm w)\bm u,\bm v\bigr)_{\Om(t)}\\
&\quad+b(t;\bm u-\bm w,\bm u,\bm v).
\end{aligned}
\end{equation}

Using \eqref{eq:ale-skew-equivalence} and \eqref{eq:DF}, together
with the boundary conditions \eqref{eq:ale-sav-traction-strong},
\eqref{eq:ale-sav-wall-strong}, \eqref{eq:ale-sav-slip-strong} and
\eqref{eq:ale-sav-young-general}, we obtain the coupled ALE--SAV weak
system
\begin{subequations}\label{eq:ale-sav-weak}
\begin{empheq}[left=\empheqlbrace]{alignat=3}
&&\;&\begin{aligned}
&(\partial_t^\circ\bm u,\bm v)_{\Om(t)}
+\frac12\bigl((\Div \bm w)\bm u,\bm v\bigr)_{\Om(t)}\\
&\quad+2\mu(\D \bm u,\D \bm v)_{\Om(t)}
+b(t;\bm u-\bm w,\bm u,\bm v)\\
&\quad
+\beta_s(P_s \bm u,P_s \bm v)_{\Ga_s(t)} +\beta_c\int_{\Ga_c(t)}(\bm u\cdot\bm\mu_s)(\bm v\cdot\bm\mu_s)\dd\ell\\
&\quad-(\pi,\Div \bm v)_{\Om(t)}+r\,\mathcal G(\bm X)[\bm v]=0
\end{aligned}\label{eq:ale-sav-momentum-continuous}&&\\
&&\;&(\Div \bm u,q)_{\Om(t)}=0.\label{eq:ale-sav-div-weak}&&\\[.8ex]
&&\;&\dot r=\tfrac12\mathcal G(\bm X)[\bm u],\label{eq:ale-sav-scalar-weak}&&
\end{empheq}
\end{subequations}
where the velocity $\bm u$ and test function $\bm v$ belong to
$H^1(\Om(t))^d$, have zero normal trace on $\Ga_s(t)$, and have
free-surface traces in $H^1(\Ga_F(t))^d$.
The pressure $\pi$ and test function $q$ belong to $L^2(\Om(t))$.

\smallskip
Define the SAV energy by
\begin{equation}\label{eq:sav-energy-continuous}
\E_{\rm SAV}(t)=\frac12\norm{\bm u}_{L^2(\Om(t))}^2+r(t)^2-C_0.
\end{equation}
For these initial data, \eqref{eq:sav-continuous-compatibility}
implies $\E_{\rm SAV}(t)=\E(t)$.

\begin{lemma}[Structural properties]\label{prop:sav-energy-continuous}
For sufficiently regular solutions of \eqref{eq:ale-sav-weak}, the
following identities hold.
\smallskip

\begin{enumerate}
\renewcommand{\labelenumi}{\textup{(\roman{enumi})}}
\item \text{\upshape Area and Volume conservation.}
\begin{equation}\label{eq:continuous-volume-invariant}
\V(\bm X(t))=\V(\bm X_0).
\end{equation}
\item \text{\upshape Energy dissipation.}
\begin{equation}\label{eq:continuous-sav-energy}
\begin{aligned}
\frac{\dd}{\dd t}\E_{\rm SAV}(t)
&=-2\mu\norm{\D\bm u}_{L^2(\Om(t))}^2
-\beta_s\norm{P_s\bm u}_{L^2(\Ga_s(t))}^2\\
&\quad-\beta_c\int_{\Ga_c(t)}|\bm u\cdot\bm\mu_s|^2\dd\ell.
\end{aligned}
\end{equation}
\end{enumerate}
\end{lemma}
\begin{proof}
\textup{(i)} The material motion gives $V_F=\bm u\cdot\bm n_F$.
With wall impermeability, \eqref{eq:continuous-volume-conservation}
expresses the volume derivative as the integral of $\Div\bm u$.
This integral vanishes by \eqref{eq:ale-sav-div-weak} with $q=1$.
\smallskip

\textup{(ii)} Test \eqref{eq:ale-sav-momentum-continuous}
with $\bm v=\bm u$ and use \eqref{eq:ale-sav-div-weak} with
$q=\pi$. The pressure term vanishes, and
$b(t;\bm u-\bm w,\bm u,\bm u)=0$ by
\eqref{eq:skew-convection-continuous}. Reynolds' transport formula gives
\begin{equation}\label{eq:ale-kinetic-continuous}
\frac{\dd}{\dd t}\frac12\norm{\bm u}_{L^2(\Om(t))}^2 = (\partial_t^\circ\bm u,\bm u)_{\Om(t)}
+\frac12\bigl((\Div\bm w)\bm u,\bm u\bigr)_{\Om(t)}.
\end{equation}
Equation~\eqref{eq:ale-sav-scalar-weak} gives
\begin{equation}
r\,\mathcal G(\bm X)[\bm u]=2r\dot r=\frac{\dd}{\dd t}r^2 .
\end{equation}
Substituting the last two identities into
\eqref{eq:ale-sav-momentum-continuous} with $\bm v=\bm u$
and using the energy definition \eqref{eq:sav-energy-continuous}
gives \eqref{eq:continuous-sav-energy}. The proof is complete.
\end{proof}

\bigskip
\section[Structure-preserving maps on evolving finite element spaces]{Structure-preserving maps on evolving finite element spaces}\label{sec:space}
In this section, we construct an $L^2$-isometric velocity transfer and a
volume-preserving boundary map on evolving finite element spaces.

\subsection{Evolving finite element spaces}
\label{sec:fe-spaces}\label{sec:isometric-transport}\label{sec:boundary-update}
Let $\mathcal T_h^0$ be an initial degree-$k$ isoparametric
simplicial mesh  \cite{LiMaQiuALE} with nodal vector $\mathbf{x}^0$.
We consider nodal deformations of $\mathcal T_h^0$ in
an admissible neighborhood, with the same element index set,
local-to-global node maps and boundary labels.
The moving nodal vector
$\mathbf{x}$ determines $\mathcal T_h(\mathbf{x})$, its domain
$\Om(\mathbf{x})$, and the boundary parts
$\Ga_F(\mathbf{x})$, $\Ga_s(\mathbf{x})$ and $\Ga_c(\mathbf{x})$.
We write $\Om_h^n:=\Om(\mathbf{x}^n)$ and
$\Ga_{j,h}^n:=\Ga_j(\mathbf{x}^n)$, $j\in\{F,s,c\}$.
For a reference simplex $\widehat K$, let $P_k(\widehat K)$ and
$P_{k-1}(\widehat K)$ denote the polynomials of total degree at most
$k$ and $k-1$, respectively.
Each element is represented by an injective map $F_K$ satisfying
\begin{equation}\label{eq:isoparametric-geometry}
K=F_K(\widehat K),\qquad F_K\in[P_k(\widehat K)]^d,\qquad
\det DF_K>0\text{ on }\widehat K.
\end{equation}
Admissible meshes satisfy \eqref{eq:isoparametric-geometry}
and are conforming, with disjoint element interiors, nondegenerate
boundary facets and contact set, positive volume, and positive
wetted-wall measure.
They also satisfy
$\overline{\Om(\mathbf{x})}\subset\{y\geq0\}$,
$\Ga_s(\mathbf{x})\subset\{y=0\}$ and
$\Ga_F(\mathbf{x})\setminus\Ga_c(\mathbf{x})\subset\{y>0\}$.
For the fluid variables, define the 
$P_k$--$P_{k-1}$ Taylor--Hood pair with $k\geq2$ by
\begin{align}
V_h(\mathbf{x})
&=\{\bm v_h\in C^0(\overline{\Om(\mathbf{x})})^d:
\bm v_h\cdot\bm n_s=0\text{ on }\Ga_s(\mathbf{x}),\nonumber\\
&\qquad\bm v_h\circ F_K\in[P_k(\widehat K)]^d
\text{ for every }K\in\mathcal T_h(\mathbf{x})\},
\label{eq:velocity-space}\\
Q_h(\mathbf{x})
&=\{q_h\in C^0(\overline{\Om(\mathbf{x})}):
q_h\circ F_K\in P_{k-1}(\widehat K)
\text{ for every }K\in\mathcal T_h(\mathbf{x})\}.
\label{eq:shifted-pressure-space}
\end{align}
The homogeneous test space for the ALE velocity is
\begin{equation}\label{eq:mdr-test-space}
\mathring V_h(\mathbf{x})
=\{\bm z_h\in V_h(\mathbf{x}):
\bm z_h=0\text{ on }\Ga_F(\mathbf{x})\}.
\end{equation}
For $k=1$, the MINI pair uses bubble-enriched continuous $P_1$ velocity
and $P_1$ pressure, with the same wall condition.
Geometry and mesh velocity are piecewise linear.

For two admissible meshes $\mathbf{x}_s,\mathbf{x}_t$
in this neighborhood, the discrete ALE map
$\Phi_{s\to t}:\Om(\mathbf{x}_s)\to\Om(\mathbf{x}_t)$ maps
corresponding geometry nodes.
For bases $\{\bm\phi_{h,j}[\mathbf{x}]\}_{j=1}^m$ of
$V_h(\mathbf{x})$ induced by the same reference shape functions and
degree-of-freedom ordering, the transport relations are
\begin{equation}\label{eq:ale-basis-transport}
\begin{aligned}
\bm\phi_{h,j}[\mathbf{x}_t]\circ\Phi_{s\to t}
&=\bm\phi_{h,j}[\mathbf{x}_s].
\end{aligned}
\end{equation}
Define the ordinary ALE transport
$\mathcal P_{s\to t}:V_h(\mathbf{x}_s)\to V_h(\mathbf{x}_t)$ by
\begin{equation}\label{eq:operator-ale-composition}
\mathcal P_{s\to t}\bm v_h
=\bm v_h\circ\Phi_{s\to t}^{-1}.
\end{equation}
By \eqref{eq:operator-ale-composition}, $\mathcal P$ is a linear
bijection, and \eqref{eq:ale-basis-transport} gives coefficient
preservation in the transported bases. It does not in general preserve
the physical $L^2$ inner product.

\subsection[Isometric transport between meshes]{Isometric transport between meshes}\label{sec:structure-operators}

For $\bm u_h\in V_h(\mathbf{x}_s)$, define the isometric transport
$\mathcal I_{s\to t}:V_h(\mathbf{x}_s)\to V_h(\mathbf{x}_t)$ by
\begin{equation}\label{eq:operator-isometric-transport}
\mathcal I_{s\to t}\bm u_h
=\sum_{j=1}^m\widetilde U_{s,j}\bm\psi_{h,j}[\mathbf{x}_t] \quad \text{with} \quad \widetilde U_{s,j}:=(\bm u_h,\bm\psi_{h,j}[\mathbf{x}_s])_{\Omega(\mathbf{x}_s)},
\end{equation}
where $\{\bm\psi_{h,j}[\mathbf{x}]\}_{j=1}^m$ is an ordered
$L^2(\Om(\mathbf{x}))$-orthonormal basis of $V_h(\mathbf{x})$,
chosen smoothly in $\mathbf{x}$ after ALE pullback.
For $\bm u_h,\bm v_h\in V_h(\mathbf{x}_s)$, orthonormality gives
\begin{equation}\label{eq:operator-transport-inner-product}
(\mathcal I_{s\to t}\bm u_h,\mathcal I_{s\to t}\bm v_h)_{\Omega(\mathbf{x}_t)}
=(\bm u_h,\bm v_h)_{\Omega(\mathbf{x}_s)}.
\end{equation}
Thus $\mathcal I$ preserves kinetic energy.
The fixed basis choices imply
$\mathcal I_{s\to s}=\mathrm{Id}$,
$\mathcal I_{s\to t}^{-1}=\mathcal I_{t\to s}$ and
$\mathcal I_{t\to t'}\mathcal I_{s\to t}=\mathcal I_{s\to t'}$
for any third admissible mesh $\mathbf{x}_{t'}$ in the same neighborhood.
% The same identity, inverse and composition relations hold for $\mathcal P$.

Along an admissible mesh path $\mathbf{x}(t)$, set
$\Om_h(t)=\Om(\mathbf{x}(t))$ and
$\Phi_h(\cdot,t)=\Phi_{0\to t}:\Om_h(0)\to\Om_h(t)$.
Assume that $\Phi_h$ is uniformly bi-Lipschitz and orientation-preserving,
with $\Phi_h\in W^{1,\infty}\bigl(0,T;W^{1,\infty}(\Om_h(0))^d\bigr)$.
Let $\bm u_h(t)\in V_h(\mathbf{x}(t))$ satisfy
$\bm u_h\circ\Phi_h\in W^{1,\infty}\bigl(0,T;L^2(\Om_h(0))^d\bigr)$,
and set $\bm w_h=\partial_t\Phi_h\circ\Phi_h^{-1}$.
The ALE derivative \eqref{eq:ale-continuous} uses
$\mathcal P_{s\to t}$ to compare velocities in the space at time $t$.
Along the same mesh path, define the time derivative associated with
the isometric transport $\mathcal I$ by
\begin{equation}\label{eq:continuous-isometric-connection}
\begin{aligned}
\partial_t^{\mathcal I}\bm u_h
&:=\left.\frac{\dd}{\dd s}
\bigl(\mathcal I_{s\to t}\bm u_h(s)\bigr)\right|_{s=t}.
\end{aligned}
\end{equation}
By \eqref{eq:operator-transport-inner-product},
$\partial_t^{\mathcal I}$ satisfies the same kinetic-energy identity
as on a fixed domain
\begin{equation}\label{eq:isometric-derivative-energy}
\frac{\dd}{\dd t}\frac12\norm{\bm u_h}_{L^2(\Om_h(t))}^2
=(\partial_t^{\mathcal I}\bm u_h,\bm u_h)_{\Om_h(t)}.
\end{equation}
Related variable transformations for energy estimates
on moving domains are studied in \cite{LundquistMalanNordstrom}.
The composition
$\mathcal I_{s\to t}\mathcal P_{t\to s}$ transports a field from time
$t$ to $s$ by $\mathcal P$ and returns it by $\mathcal I$.
Its derivative measures the difference between these transport rules.
For $\bm v_h\in V_h(\mathbf{x}(t))$, define
\begin{equation}\label{eq:isometric-connection-representation}
\mathcal D_h(t)\bm v_h
:=\left.\frac{\dd}{\dd s}
\bigl(\mathcal I_{s\to t}\mathcal P_{t\to s}\bm v_h\bigr)
\right|_{s=t}.
\end{equation}
For $\bm u_h,\bm v_h\in V_h(\mathbf{x}(t))$, define the skew-symmetric bilinear form
\begin{equation}\label{eq:geometric-skew-form}
s_h(t;\bm u_h,\bm v_h)
:=\tfrac12(\mathcal D_h(t)\bm u_h,\bm v_h)_{\Om_h(t)}
-\tfrac12(\mathcal D_h(t)\bm v_h,\bm u_h)_{\Om_h(t)}.
\end{equation}
Thus $s_h(t;\bm u_h,\bm u_h)=0$.
\smallskip

The following lemma relates the time-derivative terms
in the ALE formulation to the isometric time derivative.

\begin{lemma}\label{lem:discrete-ale-inertia}
Under the mesh and regularity assumptions stated above,
the following identities hold for almost every $t\in(0,T)$.
\begin{equation}\label{eq:isometric-derivative-relation}
\partial_t^{\mathcal I}\bm u_h
=\partial_t^\circ\bm u_h+\mathcal D_h(t)\bm u_h(t).
\end{equation}
For every $\bm v_h\in V_h(\mathbf{x}(t))$,
\begin{equation}\label{eq:isometric-ale-inertia}
% \begin{aligned}
(\partial_t^\circ\bm u_h,\bm v_h)_{\Om_h(t)}
+\tfrac12\bigl((\Div\bm w_h)\bm u_h,\bm v_h\bigr)_{\Om_h(t)} =(\partial_t^{\mathcal I}\bm u_h,\bm v_h)_{\Om_h(t)}
-s_h(t;\bm u_h,\bm v_h).
% \end{aligned}
\end{equation}
\end{lemma}
\begin{proof}
Since $\mathcal P_{t\to s}\mathcal P_{s\to t}=\mathrm{Id}$,
\begin{equation}\label{eq:isometric-transport-factorization}
\mathcal I_{s\to t}\bm u_h(s)
=\bigl(\mathcal I_{s\to t}\mathcal P_{t\to s}\bigr)
\bigl(\mathcal P_{s\to t}\bm u_h(s)\bigr).
\end{equation}
At $s=t$, the factors are $\mathrm{Id}$ and $\bm u_h(t)$,
so the product rule gives \eqref{eq:isometric-derivative-relation}.

Differentiating \eqref{eq:operator-transport-inner-product} for
$\mathcal P_{t\to s}\bm u_h$ and $\mathcal P_{t\to s}\bm v_h$
at $s=t$, Reynolds' formula gives
\begin{equation}\label{eq:connection-symmetric-part}
\begin{aligned}
&(\mathcal D_h(t)\bm u_h,\bm v_h)_{\Om_h(t)}
+(\bm u_h,\mathcal D_h(t)\bm v_h)_{\Om_h(t)} =\bigl((\Div\bm w_h)\bm u_h,\bm v_h\bigr)_{\Om_h(t)}.
\end{aligned}
\end{equation}
Combining \eqref{eq:geometric-skew-form}, \eqref{eq:isometric-derivative-relation} and
\eqref{eq:connection-symmetric-part} 
gives \eqref{eq:isometric-ale-inertia}. The proof is complete.
\end{proof}

The isometric transport also yields an exact kinetic-energy balance
for a midpoint velocity update between different meshes.

\begin{lemma}[Discrete kinetic-energy balance]\label{lem:discrete-kinetic-balance}
Let $\mathbf{x}^n,\mathbf{x}^{n+1/2},\mathbf{x}^{n+1}$ be arbitrary
admissible meshes in this neighborhood.
Given $\bm u_h^n\in V_h(\mathbf{x}^n)$ and
$\bm u_h^{n+1/2}\in V_h(\mathbf{x}^{n+1/2})$, set
\begin{equation}\label{eq:three-mesh-endpoint-recovery}
\bm u_h^{n+1}=\mathcal I_{n+1/2\to n+1}
\bigl(2\bm u_h^{n+1/2}-\mathcal I_{n\to n+1/2}\bm u_h^n\bigr).
\end{equation}
Then
\begin{equation}\label{eq:function-kinetic-increment}
\frac12 \norm{\bm u_h^{n+1}}_{L^2(\Om_h^{n+1})}^2
- \frac12 \norm{\bm u_h^n}_{L^2(\Om_h^{n})}^2 = \bigl(\bm u_h^{n+1/2}-\mathcal I_{n\to n+1/2}\bm u_h^n,
2\bm u_h^{n+1/2}\bigr)_{\Om_h^{n+1/2}}.
\end{equation}
\end{lemma}
\begin{proof}
The inverse property of $\mathcal I$ and
\eqref{eq:three-mesh-endpoint-recovery} give, in
$V_h(\mathbf{x}^{n+1/2})$,
\[
\mathcal I_{n+1\to n+1/2}\bm u_h^{n+1}
=2\bm u_h^{n+1/2}-\mathcal I_{n\to n+1/2}\bm u_h^n.
\]
Use \eqref{eq:operator-transport-inner-product} to express both
endpoint norms in the midpoint space. The difference of their
squares then gives \eqref{eq:function-kinetic-increment}. The proof is complete.
\end{proof}

\begin{remark}[Matrix form of the midpoint update]
\label{rem:isometric-transport-implementation}
The update in Lemma~\ref{lem:discrete-kinetic-balance} is an arithmetic
midpoint update in mass coordinates. For $j\in\{n,n+\tfrac12,n+1\}$,
let $U^j$ be the nodal coefficients of $\bm u_h^j$ and
$\mathbf M^j$ the exactly integrated mass matrix on $\Om_h^j$.
Using upper triangular Cholesky factors with positive diagonal, set
\begin{equation}\label{eq:solution-mass-coordinates}
\begin{gathered}
\mathbf M^j=(\mathbf R^j)^T\mathbf R^j,\qquad
\widetilde U^j=\mathbf R^jU^j,\\
\|\bm u_h^j\|_{L^2(\Om_h^j)}^2
=(U^j)^T\mathbf M^jU^j=\|\widetilde U^j\|_2^2.
\end{gathered}
\end{equation}
Choose the orthonormal bases
$\bm\psi_{h,a}^j=\sum_b\bm\phi_{h,b}^j[(\mathbf R^j)^{-1}]_{ba}$.
Then \eqref{eq:three-mesh-endpoint-recovery} becomes
\begin{equation}\label{eq:solution-isometric-transfer}
\begin{gathered}
\mathcal I_{n\to n+1/2}\bm u_h^n
=\sum_a[(\mathbf R^{n+1/2})^{-1}\widetilde U^n]_a
\bm\phi_{h,a}^{n+1/2},\\
\widetilde U^{n+1/2}
=\mathbf R^{n+1/2}U^{n+1/2}
=\tfrac12(\widetilde U^{n+1}+\widetilde U^n).
\end{gathered}
\end{equation}
For $\tau>0$, the isometric time-difference term in the weak form is
represented by
\[
(\mathbf R^{n+1/2})^T
\frac{\widetilde U^{n+1}-\widetilde U^n}{\tau}
=\frac2\tau\left[
\mathbf M^{n+1/2}U^{n+1/2}
-(\mathbf R^{n+1/2})^T\widetilde U^n\right].
\]
Testing with $U^{n+1/2}$ gives the kinetic-energy increment
\[
\begin{aligned}
(\widetilde U^{n+1/2})^T
\frac{\widetilde U^{n+1}-\widetilde U^n}{\tau}
&=\frac{\|\widetilde U^{n+1}\|_2^2
-\|\widetilde U^n\|_2^2}{2\tau} =\frac{\|\bm u_h^{n+1}\|_{L^2(\Om_h^{n+1})}^2
-\|\bm u_h^n\|_{L^2(\Om_h^n)}^2}{2\tau}.
\end{aligned}
\]
\end{remark}

\subsection[Volume-preserving boundary map]{Volume-preserving boundary map}\label{sec:volume-preserving-map}
Fix an admissible mesh \(\mathbf x\) and consider finite element boundary maps \(\mathbf Y_h:\Gamma_F(\mathbf x)\to\mathbb R^d\) whose contact nodes lie on the wall. For an admissible deformed geometry with free-surface parametrization \(\mathbf X\), the volume and its first variation under a wall-tangent deformation \(\boldsymbol\xi\) are
\begin{equation}\label{eq:area}
\mathcal V(\mathbf X)=|\Omega(t)|
=\int_{\Gamma_F(t)}y\,(\mathbf n_F\cdot\mathbf e_y)\,\mathrm ds,
\qquad
D\mathcal V(\mathbf X)[\boldsymbol\xi]
=\int_{\Gamma_F(t)}\boldsymbol\xi\cdot\mathbf n_F\,\mathrm ds.
\end{equation}
Pulling this boundary integral back to fixed oriented parameter elements defines \(\mathcal V[\mathbf Y_h]\) as a polynomial in the boundary nodal coordinates. Its derivative acts on variations with zero vertical contact components. On admissible maps, \(\mathcal V[\mathbf Y_h]\) equals the enclosed volume.

We seek an auxiliary velocity for the evolution of \(\mathbf Y_h(t)\), \(t\in[t_n,t_{n+1}]\), that preserves \(\mathcal V[\mathbf Y_h]\) and agrees with the prescribed boundary velocity at \(\mathbf Y_h=\mathrm{Id}\).

% For any admissible mesh $\mathbf{x}$, we seek an auxiliary
% velocity for the evolution of finite element boundary maps
% $\bm Y_h(t):\Ga_F(\mathbf{x})\to\mathbb R^d$,
% $t\in[t_n,t_{n+1}]$, that preserves volume and agrees with the
% prescribed velocity at $\bm Y_h=\mathrm{Id}$.

% For a finite element boundary map $\bm Y_h$ with contact nodes on the wall,
% we define $\V[\bm Y_h]$ by pulling the boundary integral in
% \eqref{eq:area} back to fixed oriented parameter elements. This gives a
% polynomial in the boundary nodal coordinates. Its derivative
% $D\V[\bm Y_h]$ acts on variations with zero vertical contact components.
% On admissible maps, $\V[\bm Y_h]$ equals the enclosed volume.

% For the wall-tangent deformations of Section~\ref{sec:geometric-sav},
% the boundary representation of the volume and its first variation are
% \begin{equation}\label{eq:area}
% \V(\bm X)=|\Om(t)|=\int_{\Ga_F(t)}y\,(\bm n_F\cdot\bm e_y)\dd s,
% \qquad
% D\V(\bm X)[\bm{\xi}]=\int_{\Ga_F(t)}\bm{\xi}\cdot\bm n_F\dd s.
% \end{equation}
% \par

\begin{lemma}[Volume-preserving boundary evolution]\label{lem:auxiliary-boundary-velocity}
Let $\bm u_{F,h}$ be a finite element boundary velocity on
$\Ga_F(\mathbf{x})$ with zero vertical contact components. Assume
$D\V[\mathrm{Id}](\bm u_{F,h})=0$ and $D\V[\mathrm{Id}]\not\equiv0$.
For every finite element boundary map $\bm Y_h$ with contact nodes
on the wall, there exists an auxiliary finite element boundary velocity
$\bm w_{F,h}[\bm Y_h]$ satisfying
\begin{equation}\label{eq:geometry-velocity-requirements}
 \bm w_{F,h}[\mathrm{Id}]=\bm u_{F,h},\qquad
 D\V[\bm Y_h]\bigl(\bm w_{F,h}[\bm Y_h]\bigr)=0.
\end{equation}
Every differentiable admissible solution of
$\dot{\bm Y}_h=\bm w_{F,h}[\bm Y_h]$ on $[t_n,t_{n+1}]$
therefore preserves volume:
\begin{equation}
 \frac{\dd}{\dd t}\V[\bm Y_h(t)]=0.
\end{equation}
\end{lemma}
\begin{proof}
Since $D\V[\mathrm{Id}]\not\equiv0$, choose a variation $\bm\xi_h$
with $D\V[\mathrm{Id}](\bm\xi_h)=1$ and define
\begin{equation}\label{eq:auxiliary-boundary-velocity}
 \bm w_{F,h}[\bm Y_h]
 =D\V[\bm Y_h](\bm\xi_h)\,\bm u_{F,h}
 -D\V[\bm Y_h](\bm u_{F,h})\,\bm\xi_h.
\end{equation}
Since $D\V[\mathrm{Id}](\bm\xi_h)=1$ and
$D\V[\mathrm{Id}](\bm u_{F,h})=0$,
\eqref{eq:auxiliary-boundary-velocity} gives
$\bm w_{F,h}[\mathrm{Id}]=\bm u_{F,h}$. By linearity,
\[
\begin{aligned}
 D\V[\bm Y_h]\bigl(\bm w_{F,h}[\bm Y_h]\bigr)
 &=D\V[\bm Y_h](\bm\xi_h)\,D\V[\bm Y_h](\bm u_{F,h}) -D\V[\bm Y_h](\bm u_{F,h})\,D\V[\bm Y_h](\bm\xi_h)=0.
\end{aligned}
\]
Volume conservation follows from
\eqref{eq:geometry-velocity-requirements} and the chain rule.
\end{proof}

We discretize this evolution over a step $\tau=t_{n+1}-t_n>0$.

\begin{lemma}[Linear volume-preserving boundary scheme]\label{lem:boundary-update-construction}
Let $\bm w_{F,h}$ be given by \eqref{eq:auxiliary-boundary-velocity}
under the assumptions of Lemma~\ref{lem:auxiliary-boundary-velocity},
with $\V$ and its derivatives evaluated exactly.
Assume that the linear systems in the construction are nonsingular.
There exists a boundary update map $\Phi_F(\tau)$, computable by
linear solves, with $\bm Y_h^{n+1}=\Phi_F(\tau)[\bm Y_h^n]$.
If $\bm Y_h^n$ and $\bm Y_h^{n+1}$ are admissible, the update
preserves volume:
\begin{equation}\label{eq:boundary-map-conservation}
 \V[\bm Y_h^{n+1}]=\V[\bm Y_h^n].
\end{equation}
\end{lemma}
\begin{proof}
On the space of variations with zero vertical contact components,
use $\langle\bm\eta_h,\bm\zeta_h\rangle_N
=\sum_i\bm\eta_h(q_i)\cdot\bm\zeta_h(q_i)$ at the boundary nodes
$q_i$. Define $\bm a_h$ in this space and set $\bm\xi_h$ by
\begin{equation}\label{eq:function-normalized-direction}
 \langle\bm a_h,\bm\eta_h\rangle_N
 =D\V[\mathrm{Id}](\bm\eta_h),\qquad
 \bm\xi_h=\frac{\bm a_h}{\langle\bm a_h,\bm a_h\rangle_N}.
\end{equation}
This choice gives $D\V[\mathrm{Id}](\bm\xi_h)=1$.
The mesh $\mathbf{x}$ and the fields $\bm u_{F,h},\bm\xi_h$
remain unchanged during this time step.

In two dimensions, define $\Phi_F(\tau)[\bm Y_h^n]$ as the midpoint
solution $\bm Y_h^{n+1}$ of
\begin{equation}\label{eq:geometry-substep}
 \frac{\bm Y_h^{n+1}-\bm Y_h^n}{\tau}
 =\bm w_{F,h}\!\left[\frac{\bm Y_h^n+\bm Y_h^{n+1}}2\right].
\end{equation}
Since area is quadratic in the boundary nodal coordinates, its
directional derivatives and hence \eqref{eq:auxiliary-boundary-velocity}
are affine in $\bm Y_h$, making \eqref{eq:geometry-substep} a linear
system. Set $\overline{\bm Y}_h=(\bm Y_h^n+\bm Y_h^{n+1})/2$.
For the quadratic area, the exact midpoint identity
\cite[p.~1039]{MQR1999}, together with
\eqref{eq:geometry-substep} and \eqref{eq:geometry-velocity-requirements},
gives
\begin{equation}\label{eq:volume-substep-identity}
 \V[\bm Y_h^{n+1}]-\V[\bm Y_h^n]
 =D\V[\overline{\bm Y}_h](\bm Y_h^{n+1}-\bm Y_h^n)=\tau D\V[\overline{\bm Y}_h]
   \bigl(\bm w_{F,h}[\overline{\bm Y}_h]\bigr)=0.
\end{equation}

In three dimensions, define the skew-symmetric operator
\begin{equation}\label{eq:function-geometry-skew-operator}
 \mathcal J_h\bm\eta_h
 =\langle\bm\xi_h,\bm\eta_h\rangle_N\bm u_{F,h}
 -\langle\bm u_{F,h},\bm\eta_h\rangle_N\bm\xi_h.
\end{equation}
Let $\mathbf P_i$ denote the coordinate projections acting pointwise
on boundary finite element fields:
$\mathbf P_i\bm\eta_h
 :=(\bm\eta_h\cdot\bm e_i)\bm e_i$, $i=1,2,3$,
where $\bm e_i$ is the $i$th Cartesian unit vector.
Use the splitting \cite[pp.~412--413]{MQ2003}
\begin{equation}\label{eq:coordinate-skew-operators}
\begin{aligned}
\mathcal K_i&=\frac12(\mathbf P_1+2\mathbf P_{i+1})
\mathcal J_h(\mathbf P_1+2\mathbf P_{i+1}),\qquad i=1,2,\\
\mathcal K_3&=-(\mathbf P_2-\mathbf P_3)\mathcal J_h(\mathbf P_2-\mathbf P_3).
\end{aligned}
\end{equation}
These operators are skew-symmetric and sum to $\mathcal J_h$.
For a substep of length $\delta>0$, obtain
$\Psi_\ell(\delta)[\bm Y_h]$ by solving, for every variation
$\bm\eta_h$,
\begin{equation}\label{eq:geometry-split-substep}
\left\langle
\frac{\Psi_\ell(\delta)[\bm Y_h]-\bm Y_h}{\delta},\bm\eta_h
\right\rangle_N
=-D\V\!\left[\frac{\bm Y_h+\Psi_\ell(\delta)[\bm Y_h]}2\right]
(\mathcal K_\ell\bm\eta_h).
\end{equation}
Each $\mathcal K_\ell$ acts on two coordinate groups, so testing
with variations supported in the third group shows that this group
remains fixed. In \eqref{eq:area}, the height is linear in the vertical
group and the oriented projected area is bilinear in the horizontal
groups. Thus $\V$ is at most quadratic in the active groups. The substep is linear and
preserves $\V$ by the quadratic midpoint identity used in
\eqref{eq:volume-substep-identity} and the skew symmetry of
$\mathcal K_\ell$. The symmetric composition
\begin{equation}\label{eq:volume-split-order}
\Phi_F(\tau)=
\Psi_1(\tau/2)\circ\Psi_2(\tau/2)\circ\Psi_3(\tau)
\circ\Psi_2(\tau/2)\circ\Psi_1(\tau/2)
\end{equation}
also preserves $\V$. By \eqref{eq:auxiliary-boundary-velocity},
\[
\langle\bm w_{F,h}[\bm Y_h],\bm\eta_h\rangle_N
=-D\V[\bm Y_h](\mathcal J_h\bm\eta_h).
\]
Since $\sum_\ell\mathcal K_\ell=\mathcal J_h$, the split fields
sum to $\bm w_{F,h}$. The symmetric composition of midpoint
substeps is therefore consistent and second order for the auxiliary
evolution with the data fixed. In two dimensions, the midpoint rule
also gives second order.
\end{proof}

\section[The fully discrete CN-ALE--SAV finite element method]{The fully discrete CN-ALE--SAV finite element method}\label{sec:scheme}\label{sec:time-step}\label{sec:pressure}\label{sec:discrete-energy}
In this section, we combine the SAV formulation with the structure-preserving maps
to construct a fully discrete scheme requiring only linear solves
at each time step.
For the SAV coupling, define the discrete geometric energy by
\begin{equation}\label{eq:discrete-geometric-energy}
\E_{{\rm G},h}(\mathbf{x})
=\gamma|\Ga_F(\mathbf{x})|-\gamma\cos\theta_Y|\Ga_s(\mathbf{x})|
+g\int_{\Om(\mathbf{x})}y\dd x.
\end{equation}

For $\bm v_h\in V_h(\mathbf{x})$, the notation
$D\E_{{\rm G},h}(\mathbf{x})[\bm v_h]$ denotes the derivative under nodal
perturbations $\mathbf{x}_i+\varepsilon\bm v_h(\mathbf{x}_i)$,
with $\bm v_h$ fixed at the base geometry.
For $C_0>0$, define the functional $\mathcal G_h(\mathbf{x})$ by
\begin{equation}\label{eq:sav}
\mathcal G_h(\mathbf{x})[\bm v_h]
=\frac{D\E_{{\rm G},h}(\mathbf{x})[\bm v_h]}
{\sqrt{\E_{{\rm G},h}(\mathbf{x})+C_0}}.
\end{equation}
% The momentum and SAV updates use this same functional.

Let $t_n=n\tau$, where $\tau>0$, and choose an admissible initial
mesh $\mathbf{x}^0$. The initial velocity $\bm u_h^0\in V_h(\mathbf{x}^0)$
satisfies $(\Div\bm u_h^0,q_h)_{\Om_h^{0}}=0$ for every
$q_h\in Q_h(\mathbf{x}^0)$. Initialize
% \begin{equation}\label{eq:function-initial-state}
$r^0=\sqrt{\E_{{\rm G},h}(\mathbf{x}^0)+C_0}$.
% \end{equation}
For $\bm u_h\in V_h(\mathbf{x})$, the discrete ALE harmonic
extension $\bm w_h=\MDR_{h,\mathbf{x}}[\bm u_h]$ belongs to
$\bm u_h+\mathring V_h(\mathbf{x})$ and satisfies the finite element
form of \eqref{eq:mdr-continuous-weak}:
\begin{equation}\label{eq:mdr-discrete-weak}
(\D\bm w_h,\D\bm z_h)_{\Om(\mathbf{x})}=0
\end{equation}
for every $\bm z_h\in\mathring V_h(\mathbf{x})$.
For the MINI pair, the trial and test fields for the ALE harmonic
extension are restricted to their piecewise linear subspaces. The prescribed
free-surface trace is piecewise linear because element bubbles
vanish on element faces.
At time $t_n$, the state $(\mathbf{x}^n,\bm u_h^n,r^n)$, with
$\bm u_h^n\in V_h(\mathbf{x}^n)$, determines the discrete ALE velocity
\begin{equation}\label{eq:function-instantaneous-mdr}
\bm w_h^n=\MDR_{h,\mathbf{x}^n}[\bm u_h^n] \in\bm u_h^n +\mathring V_h(\mathbf{x}^n) 
\end{equation}
Compute $\bm w_h^0$ from \eqref{eq:function-instantaneous-mdr} before Step~1.
Each step solves the fluid--SAV system on a predicted
midpoint mesh, advances the free surface and extends its displacement
to the bulk, and recovers the endpoint velocity by isometric transport.

\algorithmstep{Step 1: Predicted geometry and velocities}\label{sec:mdr-discrete}
A hat denotes predicted data. For fixed $h$, the local estimates below
use smooth exact histories and smooth transfer operators.
On $V_h(\mathbf{x}^n)$, use the second-order midpoint extrapolation
for $n\ge1$:
\begin{equation}\label{eq:function-mesh-predictor}
\widehat{\bm w}_h^{\,n+1/2}=
\begin{cases}
\bm w_h^0,&n=0,\\
\tfrac32\bm w_h^n-\tfrac12
\mathcal P_{n-1\to n}\bm w_h^{n-1},&n\geq1.
\end{cases}
\end{equation}
Use this field to predict the midpoint and auxiliary endpoint:
\begin{equation}\label{eq:function-predicted-geometries}
\begin{aligned}
\widehat{\mathbf{x}}_i^{\,n+1/2}
=\mathbf{x}_i^n+\tfrac\tau2
\widehat{\bm w}_h^{\,n+1/2}(\mathbf{x}_i^n) \quad \text{and} \quad
\widehat{\mathbf{x}}_i^{n+1}
=\mathbf{x}_i^n+\tau
\widehat{\bm w}_h^{\,n+1/2}(\mathbf{x}_i^n).
\end{aligned}
\end{equation}
Require both predicted meshes to be admissible.
The predicted midpoint geometry has an $O(\tau^2)$ position error.
Write $\widehat\Om_h^{\,n+1/2}=\Om(\widehat{\mathbf{x}}^{\,n+1/2})$ and
$\widehat\Ga_{j,h}^{\,n+1/2}=\Ga_j(\widehat{\mathbf{x}}^{\,n+1/2})$,
$j=F,s,c$. Half-integer labels in $\mathcal P$ and $\mathcal I$
refer to this predicted midpoint mesh. The auxiliary predicted endpoint
$\widehat{\mathbf{x}}^{\,n+1}$ is used only in $\mathcal R_h(t_{n+1})$
in \eqref{eq:function-isometry-difference}.

Transport the fluid histories to the midpoint mesh and apply the same
second-order extrapolation for $n\ge1$:
\begin{equation}\label{eq:function-advector}
\begin{aligned}
\widehat{\bm u}_h^{\,n+1/2}
&=\begin{cases}
\mathcal P_{0\to1/2}\bm u_h^0,&n=0,\\
\tfrac32\mathcal P_{n\to n+1/2}\bm u_h^n
-\tfrac12\mathcal P_{n-1\to n+1/2}\bm u_h^{n-1},&n\geq1,
\end{cases}\\
\widehat{\boldsymbol\beta}_h^{\,n+1/2}
&=\widehat{\bm u}_h^{\,n+1/2}
-\mathcal P_{n\to n+1/2}\widehat{\bm w}_h^{\,n+1/2}.
\end{aligned}
\end{equation}
These fields lie in $V_h(\widehat{\mathbf{x}}^{\,n+1/2})$.
By the trace condition \eqref{eq:function-instantaneous-mdr} and the
composition property of $\mathcal P$, these extrapolations give
$\widehat{\boldsymbol\beta}_h^{\,n+1/2}=0$ on
$\widehat\Ga_{F,h}^{\,n+1/2}$.

To form the skew-symmetric correction in \eqref{eq:isometric-ale-inertia},
use the predicted affine path
$\widehat{\mathbf x}_i(s)=\mathbf x_i^n+(s-t_n)
\widehat{\bm w}_h^{\,n+1/2}(\mathbf x_i^n)$,
$s\in[t_n,t_{n+1}]$. For admissible $\widehat{\mathbf x}(s)$, define
on the fixed space $V_h(\widehat{\mathbf{x}}^{\,n+1/2})$
\begin{align}
\mathcal D_h^{n+1/2}
:=\frac{\mathcal R_h(t_{n+1})-\mathcal R_h(t_n)}{\tau} \quad\text{with}\quad\mathcal R_h(s)
=\mathcal I_{s\to t_{n+1/2}}\mathcal P_{t_{n+1/2}\to s}.
\label{eq:function-isometry-difference}
\end{align}
Along a smooth admissible predicted path, the centered difference
approximates the derivative in \eqref{eq:isometric-connection-representation}
at $t_{n+1/2}$ with error $O(\tau^2)$.
For $\bm u_h,\bm v_h\in V_h(\widehat{\mathbf x}^{\,n+1/2})$,
define the skew-symmetric form used in Step~2 by
\begin{equation}\label{eq:function-geometric-skew-form}
\begin{aligned}
s_h^{n+1/2}(\bm u_h,\bm v_h)
:={}&\tfrac12\bigl(\mathcal D_h^{n+1/2}\bm u_h,\bm v_h\bigr)_{\widehat\Om_h^{\,n+1/2}}
-\tfrac12\bigl(\mathcal D_h^{n+1/2}\bm v_h,\bm u_h\bigr)_{\widehat\Om_h^{\,n+1/2}}.
\end{aligned}
\end{equation}
For the coefficients in Step~2, write
$b_h(\mathbf x;\cdot,\cdot,\cdot)$ for the finite element
restriction of \eqref{eq:skew-convection-continuous} to $\Om(\mathbf x)$,
and set $b_h^{n+1/2}:=b_h(\widehat{\mathbf x}^{\,n+1/2};\cdot,\cdot,\cdot)$.
With $\E_{{\rm G},h}(\widehat{\mathbf{x}}^{\,n+1/2})+C_0>0$,
evaluate $\mathcal G_h(\widehat{\mathbf{x}}^{\,n+1/2})$ from
\eqref{eq:sav}.

\algorithmstep[Step 2: Linearized Crank--Nicolson fluid--SAV solve]{Step 2: Linearized Crank--Nicolson fluid--SAV solve}
The linearized Crank--Nicolson step uses isometric transport for the
velocity midpoint, the advecting field \eqref{eq:function-advector}
extrapolated to second order for $n\ge1$, and the SAV functional
evaluated on the predicted midpoint mesh.

Find the midpoint velocity
$\bm u_h^{n+1/2}\in V_h(\widehat{\mathbf{x}}^{\,n+1/2})$,
the shifted stage pressure
$\pi_h^{n+1/2}\in Q_h(\widehat{\mathbf{x}}^{\,n+1/2})$, and
$r^{n+1}\in\R$, with
\begin{equation}\label{eq:function-midpoint-values}
r^{n+1/2}=\tfrac12(r^{n+1}+r^n).
\end{equation}

For every $\bm v_h\in V_h(\widehat{\mathbf{x}}^{\,n+1/2})$ and
$q_h\in Q_h(\widehat{\mathbf{x}}^{\,n+1/2})$, require
\begin{subequations}\label{eq:variational-step}
\begin{empheq}[left=\empheqlbrace]{align}
&\begin{aligned}
&\frac2\tau\bigl(\bm u_h^{n+1/2}
-\mathcal I_{n\to n+1/2}\bm u_h^n,\bm v_h\bigr)_{\widehat\Om_h^{\,n+1/2}} -s_h^{n+1/2}(\bm u_h^{n+1/2},\bm v_h)\\[.4ex]
&\qquad+b_h^{n+1/2}\bigl(
\widehat{\boldsymbol\beta}_h^{\,n+1/2},\bm u_h^{n+1/2},\bm v_h\bigr)+2\mu(\D\bm u_h^{n+1/2},\D\bm v_h)_{\widehat\Om_h^{\,n+1/2}}\\
&\qquad
+\beta_s(P_s\bm u_h^{n+1/2},P_s\bm v_h)_{\widehat\Ga_{s,h}^{\,n+1/2}} +\beta_c(\bm u_h^{n+1/2}\cdot\bm\mu_s,
\bm v_h\cdot\bm\mu_s)_{\widehat\Ga_{c,h}^{\,n+1/2}}\\
&\qquad-(\pi_h^{n+1/2},\Div\bm v_h)_{\widehat\Om_h^{\,n+1/2}}
+r^{n+1/2}\mathcal G_h(\widehat{\mathbf{x}}^{\,n+1/2})[\bm v_h]=0,
\end{aligned}\label{eq:variational-momentum}\\
&(\Div\bm u_h^{n+1/2},q_h)_{\widehat\Om_h^{\,n+1/2}}=0,
\label{eq:variational-divergence}\\
&r^{n+1}-r^n
=\frac\tau2\mathcal G_h(\widehat{\mathbf{x}}^{\,n+1/2})[\bm u_h^{n+1/2}].
\label{eq:variational-sav}
\end{empheq}
\end{subequations}
The factor $2/\tau$ follows from the midpoint difference
in Remark~\ref{rem:isometric-transport-implementation},
applied on the predicted midpoint mesh.

\algorithmstep{Step 3: Volume-preserving mesh update}
For the construction in Section~\ref{sec:volume-preserving-map}, take
$\mathbf{x}=\widehat{\mathbf{x}}^{\,n+1/2}$.
The boundary map
$\bm Y_h:\widehat\Ga_{F,h}^{\,n+1/2}\to\R^d$
maps the predicted midpoint free surface to a candidate free surface.
Here $\mathrm{Id}$ represents the predicted midpoint geometry,
and $\V[\bm Y_h]$ is the enclosed volume for admissible deformations.
The step starts from $\bm Y_h^n$, with prescribed velocity
$\bm u_{F,h}$:
\begin{equation}\label{eq:boundary-old-map}
 \bm Y_h^n
 =\Phi_{n+1/2\to n}|_{\widehat\Ga_{F,h}^{\,n+1/2}},\qquad
 \bm u_{F,h}
 =\bm u_h^{n+1/2}|_{\widehat\Ga_{F,h}^{\,n+1/2}}.
\end{equation}
By \eqref{eq:area}, $D\V[\mathrm{Id}](\bm u_{F,h})$ is the
free-surface flux. It vanishes by \eqref{eq:variational-divergence}
with $q_h=1$ and wall impermeability.
\[
D\V[\mathrm{Id}](y\bm e_y)=|\widehat\Om_h^{\,n+1/2}|>0.
\]
Choose $\bm\xi_h$ by \eqref{eq:function-normalized-direction} and
compute the new boundary over the full step $\tau$ by
Lemma~\ref{lem:boundary-update-construction}:
\begin{equation}\label{eq:boundary-map-update}
 \bm Y_h^{n+1}=\Phi_F(\tau)[\bm Y_h^n].
\end{equation}

Use the boundary displacement divided by $\tau$ as the free-surface
data to solve \eqref{eq:mdr-discrete-weak} on
$\widehat\Om_h^{\,n+1/2}$, and update all geometry nodes with the
resulting mesh velocity:
\begin{equation}\label{eq:corrected-mesh-update}
\begin{aligned}
\left.\overline{\bm w}_h^{\,n+1/2}\right|_{\widehat\Ga_{F,h}^{\,n+1/2}}
&=\frac{\bm Y_h^{n+1}-\bm Y_h^n}{\tau},\\
\mathbf x_i^{n+1}
&=\mathbf x_i^n+\tau\overline{\bm w}_h^{\,n+1/2}
(\widehat{\mathbf x}_i^{\,n+1/2}),\qquad \text{for each geometry node }i.
\end{aligned}
\end{equation}

\algorithmstep{Step 4: Endpoint velocity update}
% Require the updated mesh $\mathbf{x}^{n+1}$ from
% \eqref{eq:corrected-mesh-update} to be admissible. 
After updating the geometry nodes by
\eqref{eq:corrected-mesh-update}, require the resulting mesh
$\mathbf{x}^{n+1}$ to be admissible. Retain $r^{n+1}$ from
\eqref{eq:variational-sav} and compute the endpoint velocity
$\bm u_h^{n+1}\in V_h(\mathbf{x}^{n+1})$ by
% Retain $r^{n+1}$
% from \eqref{eq:variational-sav} and complete the CN step by setting
\begin{equation}\label{eq:function-accepted-state}
\bm u_h^{n+1}
=\mathcal I_{n+1/2\to n+1}
\bigl(2\bm u_h^{n+1/2}-\mathcal I_{n\to n+1/2}\bm u_h^n\bigr).
\end{equation}
The kinetic-energy balance is \eqref{eq:function-kinetic-increment}.
Recompute $\bm w_h^{n+1}$ by \eqref{eq:function-instantaneous-mdr} and
retain $(\mathbf{x}^j,\bm u_h^j,\bm w_h^j)$, $j=n,n+1$, for Step~1.

\section[Discrete structure-preserving properties]{Discrete structure-preserving properties}
\label{sec:proofs}\label{sec:energy-identity}\label{sec:area} 
In this section, we establish a discrete SAV energy-dissipation law and exact area or volume conservation for the fully discrete scheme.
The argument for
\eqref{eq:geometric-energy-lower-bound} also applies to the discrete
geometry, giving $\E_{{\rm G},h}(\mathbf x)\geq0$ and a SAV denominator
bounded below by $\sqrt{C_0}>0$.

Define the discrete physical and SAV energies by
% \begin{equation}
\begin{align}\label{eq:Emod}
 \E_h^n
 &=\tfrac12\norm{\bm u_h^n}_{L^2(\Om_h^n)}^2
 +\E_{{\rm G},h}(\mathbf x^n),\\
 \E_{{\rm SAV},h}^n
 &=\tfrac12\norm{\bm u_h^n}_{L^2(\Om_h^n)}^2+(r^n)^2-C_0.
\end{align}
% \end{equation}

% \begin{samepage}

\begin{theorem}[Discrete SAV energy dissipation and volume conservation]
\label{thm:energy}\label{thm:area}\label{lem:boundary-volume-map}
Assume that the current, predicted and final meshes are admissible
and the geometric substeps are nonsingular. Then each well-defined
time step has the following properties.
\begin{enumerate}
\renewcommand{\labelenumi}{\textup{(\roman{enumi})}}
\item \textup{SAV energy dissipation.}
\begin{equation}\label{eq:energy}
\E_{{\rm SAV},h}^{n+1}\leq\E_{{\rm SAV},h}^n
\end{equation}
\item \textup{Area or volume conservation.}
\begin{equation}\label{eq:discrete-volume-conservation}
|\Om_h^{n+1}|=|\Om_h^n|.
\end{equation}
\end{enumerate}
\end{theorem}
% \end{samepage}
\begin{proof}
(i) \textit{SAV energy dissipation.}
Take $\bm v_h=\bm u_h^{n+1/2}$ in
\eqref{eq:variational-momentum} and
$q_h=\pi_h^{n+1/2}$ in \eqref{eq:variational-divergence}.
The divergence constraint and skew symmetry give
\[
\begin{aligned}
 (\pi_h^{n+1/2},\Div\bm u_h^{n+1/2})
 _{\widehat\Om_h^{\,n+1/2}}&=0,\\
 s_h^{n+1/2}(\bm u_h^{n+1/2},\bm u_h^{n+1/2})&=0,\\
 b_h^{n+1/2}\bigl(\widehat{\boldsymbol\beta}_h^{\,n+1/2},
 \bm u_h^{n+1/2},\bm u_h^{n+1/2}\bigr)&=0.
\end{aligned}
\]
By \eqref{eq:function-kinetic-increment},
the endpoint update \eqref{eq:function-accepted-state}
and the $L^2$-isometry of $\mathcal I$, we obtain
\begin{equation}\label{eq:kinetic-midpoint-identity}
\begin{aligned}
 &\frac{2}{\tau}\bigl(
 \bm u_h^{n+1/2}-\mathcal I_{n\to n+1/2}\bm u_h^n,
 \bm u_h^{n+1/2}\bigr)_{\widehat\Om_h^{\,n+1/2}}
%  &\quad=
%  \frac{
%  \norm{2\bm u_h^{n+1/2}
%  -\mathcal I_{n\to n+1/2}\bm u_h^n}
%  _{L^2(\widehat\Om_h^{\,n+1/2})}^2
%  -\norm{\mathcal I_{n\to n+1/2}\bm u_h^n}
%  _{L^2(\widehat\Om_h^{\,n+1/2})}^2}{2\tau}\\
=
 \frac{\norm{\bm u_h^{n+1}}_{L^2(\Om_h^{n+1})}^2
 -\norm{\bm u_h^n}_{L^2(\Om_h^n)}^2}{2\tau}.
\end{aligned}
\end{equation}
Similarly, \eqref{eq:variational-sav} and
\eqref{eq:function-midpoint-values} yield
\begin{equation}\label{eq:chain}
\begin{aligned}
 r^{n+1/2}\mathcal G_h(\widehat{\mathbf x}^{\,n+1/2})
 [\bm u_h^{n+1/2}]
 &=\frac{2r^{n+1/2}(r^{n+1}-r^n)}{\tau}=\frac{(r^{n+1})^2-(r^n)^2}{\tau}.
\end{aligned}
\end{equation}
Substituting \eqref{eq:kinetic-midpoint-identity}-\eqref{eq:chain} into the tested momentum equation and
multiplying by $\tau$ gives
\begin{equation}\label{eq:discrete-sav-energy-balance}
\begin{aligned}
 &\E_{{\rm SAV},h}^{n+1}-\E_{{\rm SAV},h}^n
 +2 \tau \mu\norm{\D\bm u_h^{n+1/2}}
 _{L^2(\widehat\Om_h^{\,n+1/2})}^2\\
 &\qquad
 +\tau \beta_s\norm{P_s\bm u_h^{n+1/2}}
 _{L^2(\widehat\Ga_{s,h}^{\,n+1/2})}^2
 +\tau \beta_c\norm{\bm u_h^{n+1/2}\cdot\bm\mu_s}
 _{L^2(\widehat\Ga_{c,h}^{\,n+1/2})}^2=0.
\end{aligned}
\end{equation}
The nonnegativity of the dissipation terms gives \eqref{eq:energy}.

(ii) \textit{Exact area or volume conservation.}
By \eqref{eq:boundary-old-map} and
\eqref{eq:corrected-mesh-update}, every free-surface geometry node
satisfies
\[
\begin{aligned}
\bm Y_h^n\bigl(\widehat{\mathbf x}_i^{\,n+1/2}\bigr)
&=\mathbf x_i^n,\\
\bm Y_h^{n+1}\bigl(\widehat{\mathbf x}_i^{\,n+1/2}\bigr)
&=\mathbf x_i^n
+\tau\overline{\bm w}_h^{\,n+1/2}
 \bigl(\widehat{\mathbf x}_i^{\,n+1/2}\bigr)
=\mathbf x_i^{n+1}.
\end{aligned}
\]
Since the boundary maps use the same isoparametric representation,
these nodal identities give
\[
\bm Y_h^n\bigl(\widehat\Ga_{F,h}^{\,n+1/2}\bigr)
=\Ga_{F,h}^n,
\qquad
\bm Y_h^{n+1}\bigl(\widehat\Ga_{F,h}^{\,n+1/2}\bigr)
=\Ga_{F,h}^{n+1}.
\]
On the admissible endpoint meshes, $\V$ therefore equals the
corresponding enclosed volume. Using
\eqref{eq:boundary-map-update} and
\eqref{eq:boundary-map-conservation}, we obtain
\begin{equation}\label{eq:corrected-volume-identity}
\begin{aligned}
|\Om_h^{n+1}|
&=\V[\bm Y_h^{n+1}]
=\V\bigl[\Phi_F(\tau)[\bm Y_h^n]\bigr]=\V[\bm Y_h^n]
=|\Om_h^n|.
\end{aligned}
\end{equation}
The proof is complete.
\end{proof}

\section{Numerical experiments}\label{sec:diagnostic}\label{sec:numerics-2d}
\suppressfloats[t]
We present numerical experiments for the proposed ALE-SAV FEM in two
and three spatial dimensions. Unless otherwise stated, we use the SAV
energy shift $C_0=1$, the surface tension coefficient $\gamma=1$,
and the mass factorization in \eqref{eq:solution-mass-coordinates}.
All computations are performed with NGSolve (https://ngsolve.org).
% In the two- and three-dimensional convergence tests, $g=0$, so the
% stage pressure is $p_h=\pi_h$.

\subsection{Two-dimensional droplet problem}\label{sec:ghl-energy-2d}
In the convergence tests, the initial droplet occupies the upper
half-ellipse
\begin{equation}\label{eq:convergence-initial-domain}
\Om(0)=\left\{(x,y):\frac{x^2}{a^2}+\frac{y^2}{b^2}<1,\ y>0\right\},
\qquad \bm u_0=\bm 0,
\end{equation}
with horizontal and vertical semiaxes $a=1.5$ and $b=1$,
respectively. We set $\mu=0.1$, $g=\beta_s=\beta_c=0$,
$\theta_Y=\pi/2$, and $T=0.1$.

\begin{table}[H]
\centering
\small
\caption{Finite element spaces for the two-dimensional convergence tests.}
\label{tab:convergence-elements}
\begin{tabular}{@{}ccccc@{}}
\toprule
$k$ & \shortstack{Fluid velocity $\bm u_h$}
& \shortstack{Pressure $p_h$}
& \shortstack{Geometry/ALE map $\Phi_h$}
& \shortstack{Mesh velocity $\bm w_h$}\\
\midrule
1 & $[P_1\oplus B_3]^2$ (MINI) & $P_1$ & $[P_1]^2$ & $[P_1]^2$\\
2 & $[P_2]^2$ & $P_1$ & $[P_2]^2$ & $[P_2]^2$\\
3 & $[P_3]^2$ & $P_2$ & $[P_3]^2$ & $[P_3]^2$\\
\bottomrule
\end{tabular}
\end{table}

\begin{figure}[H]
\centering
\includegraphics[width=.43\linewidth]{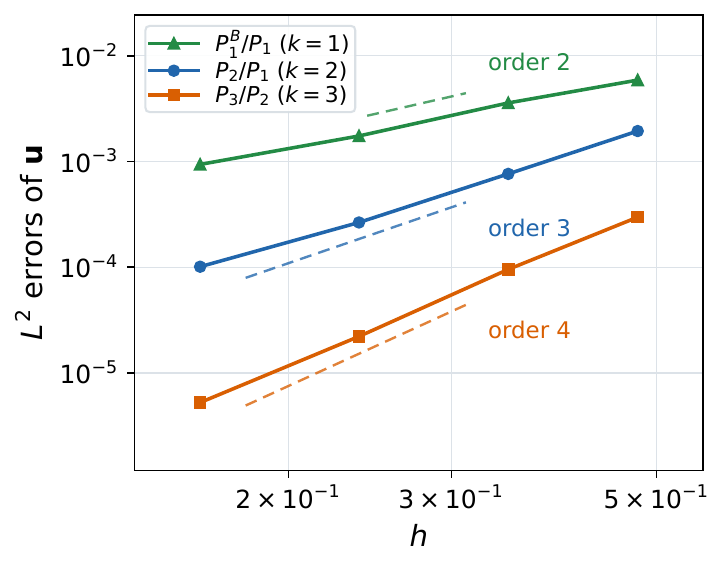}
\includegraphics[width=.43\linewidth]{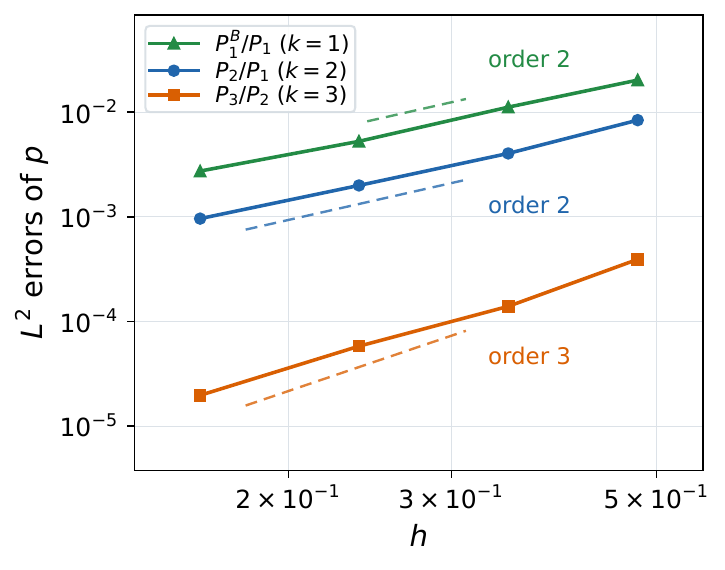}\\[-5pt]
\includegraphics[width=.43\linewidth]{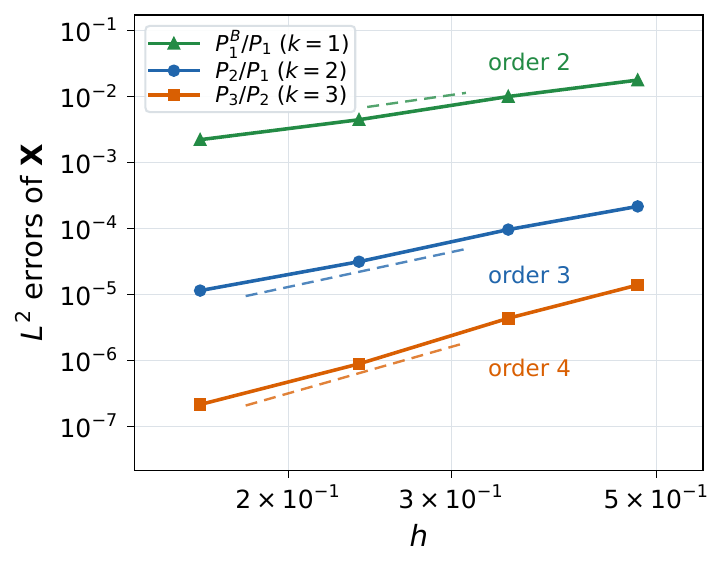}
\includegraphics[width=.43\linewidth]{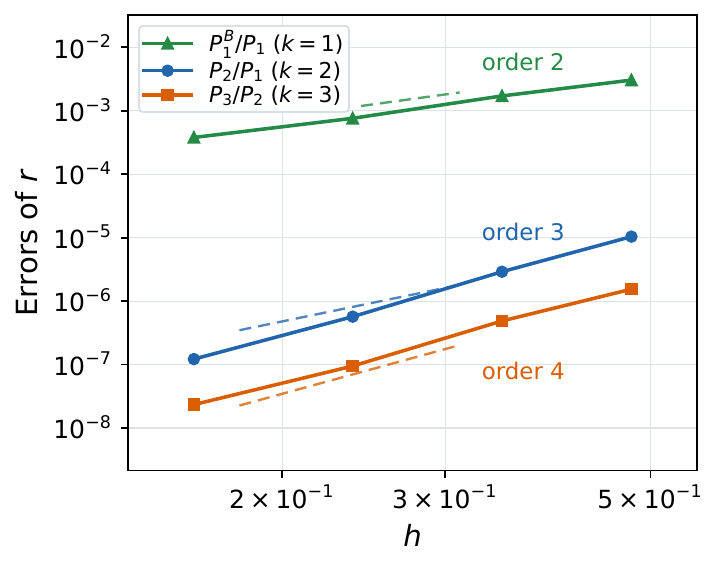}\\[-5pt]
\includegraphics[width=.43\linewidth]{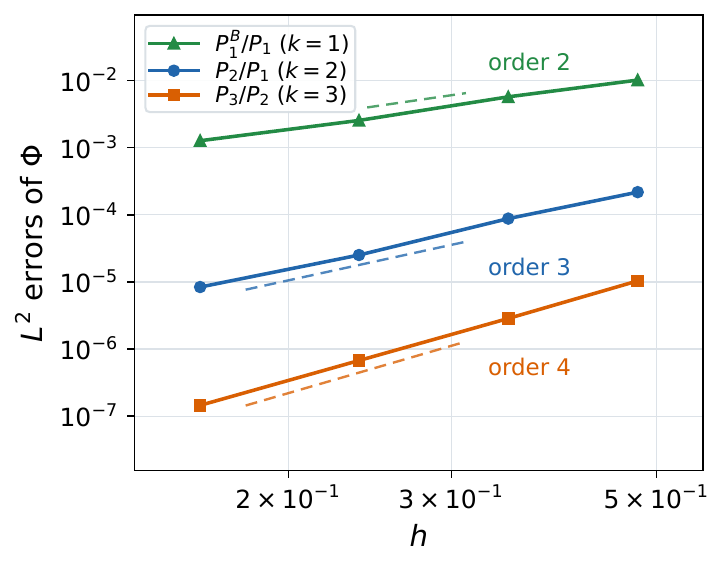}
\includegraphics[width=.43\linewidth]{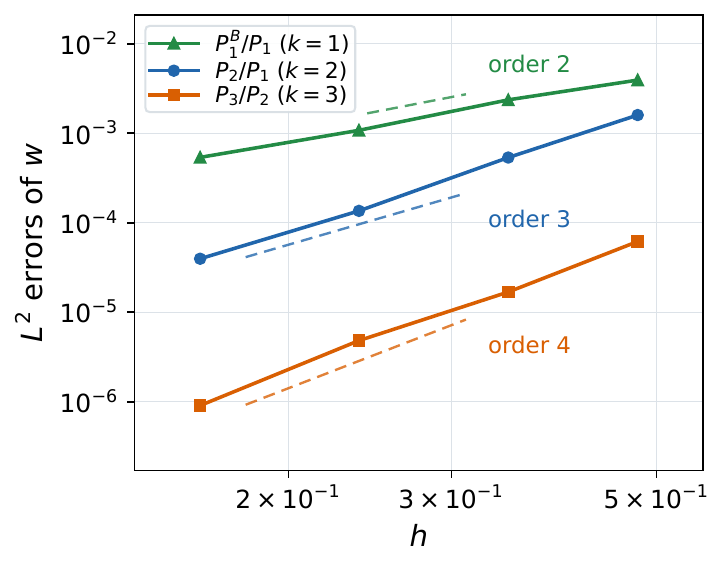}\\[-5pt]
\caption{Spatial discretization errors of the two-dimensional numerical solutions.}
\label{fig:spatial-2d}
% \end{figure}
\vspace{10pt}
% \begin{figure}[t]
\centering
\includegraphics[width=.43\linewidth]{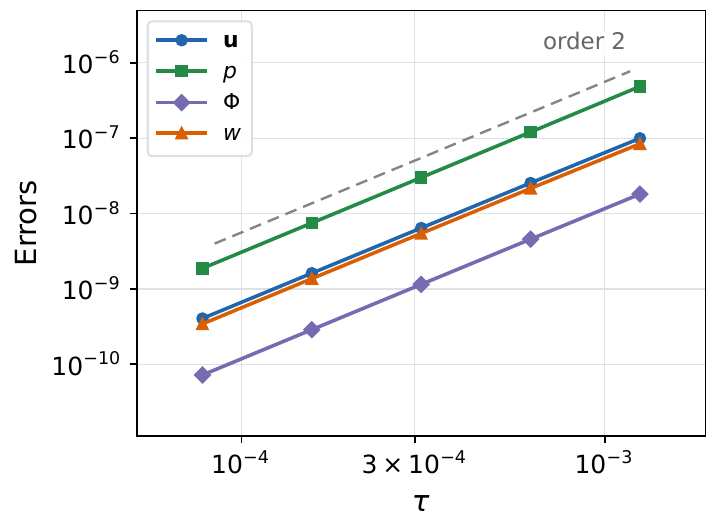}
\includegraphics[width=.43\linewidth]{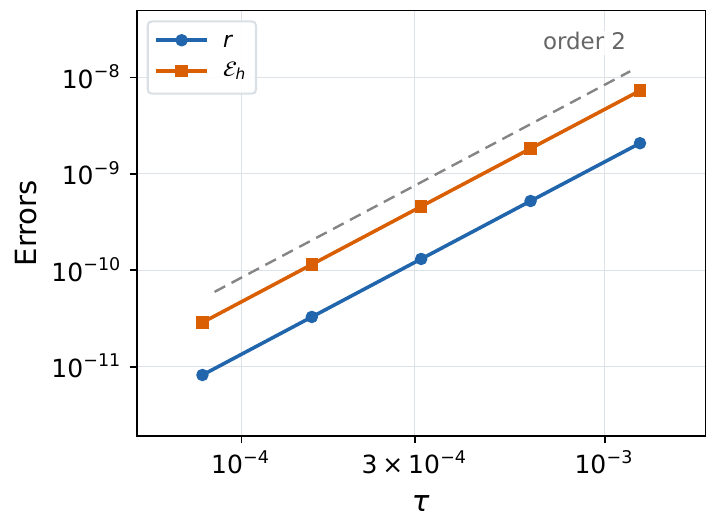}\\[.5ex]
\vspace{-5pt}
\caption{Temporal discretization errors of the two-dimensional numerical solutions.}
\label{fig:temporal-2d}
\end{figure}

Figure~\ref{fig:spatial-2d} shows the spatial discretization errors
at the final time $T=0.1$. We use the finite element spaces listed
in Table~\ref{tab:convergence-elements} for $k=1,2,3$ and fix the
time step size at $\tau=3\times10^{-4}$ so that the temporal errors
are negligible relative to the spatial errors. The reference
solution is computed using degree-$4$ finite elements on a mesh
with $h_{\rm ref}=0.1$ and the same time step size. Under mesh
refinement, the observed convergence rates for
$\bm u,\bm X,\Phi,\bm w$, and $r$ are close to $2,3,4$ for
$k=1,2,3$, respectively. 
% The corresponding pressure convergence
% rates are close to $2,2,3$. 

Figure~\ref{fig:temporal-2d} shows that the
errors in all quantities decrease at second order as the time
step is refined. To examine temporal convergence
separately, we use $P_3/P_2$ elements and fix the mesh size at
$h=0.15$ so that the spatial errors are negligible relative to
the temporal errors.

% The spatial discretization errors of the numerical solutions at $T = 0.1$ are presented in Figure~\ref{fig:spatial-2d}, where we have used the finite elements of degree $k=1,2,3$ listed in Table~\ref{tab:convergence-elements}, with a sufficiently small time step sizes $\tau = 3 \times 10^{-4}$ so that the temporal discretization errors are negligibly small compared to the spatial errors. 
% The reference solution is computed with finite elements of degree $k=4$ with $h_{\rm ref}=0.1$, and the same time step size. 
% The numerical results indicate that the convergence rates for
% $\bm u,\bm X,\Phi,\bm w, r$ are close to $2,3,4$ for
% $k=1,2,3$, respectively. The corresponding rates are approximately
% $2,2,3$ for pressure.

We next compare the numerical solutions 
for different gravitational constants $g$ and contact angles $\theta_Y$. We take the initial droplet to be a semicircle by
setting $a=b=1$ in \eqref{eq:convergence-initial-domain}, with an
initial contact angle close to $90^\circ$. We use $P_2/P_1$ elements
with $h=0.1$ and set $\beta_c=0.01$, $\beta_s=0.1$, and
$\tau=0.0025$.

\begin{figure}[H]
\centering
\includegraphics[width=.98\linewidth]{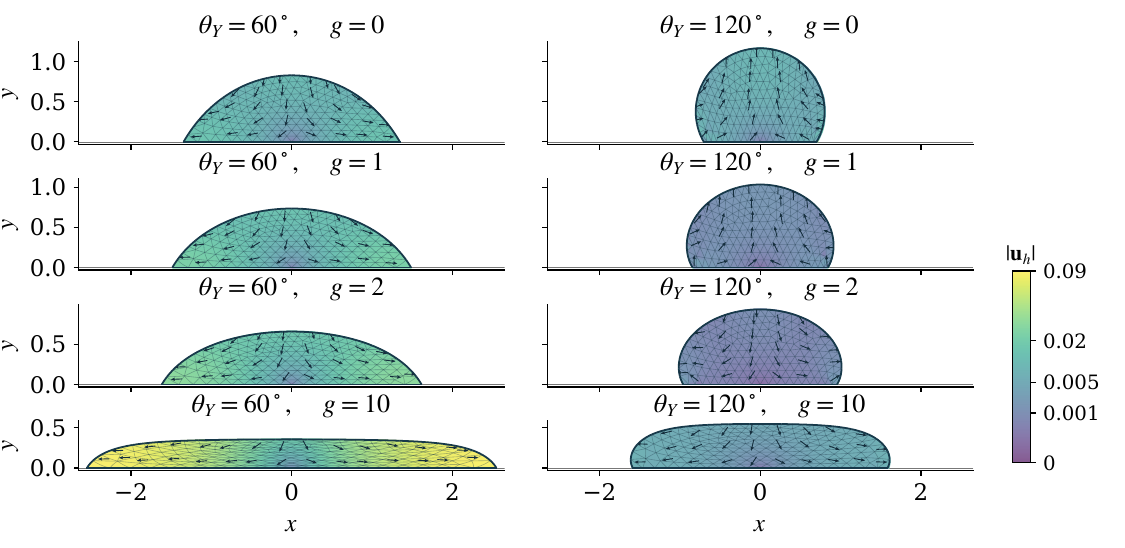}
\vspace{-5pt}
\caption{Meshes and velocities at $t = 5$ for two-dimensional wetting and
dewetting tests under different gravitational strengths.}
\label{fig:ghl-gravity-meshes}
% \end{figure}
\vspace{10pt}
% \begin{figure}[t]
\centering
\begin{subfigure}[b]{0.48\linewidth}
\centering
\includegraphics[width=\linewidth]{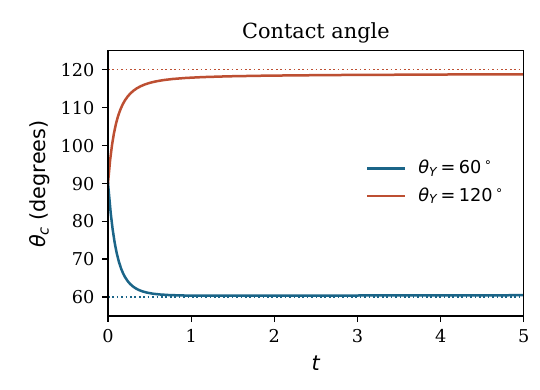}
\caption{}
\end{subfigure}
\begin{subfigure}[b]{0.48\linewidth}
\centering
\includegraphics[width=\linewidth]{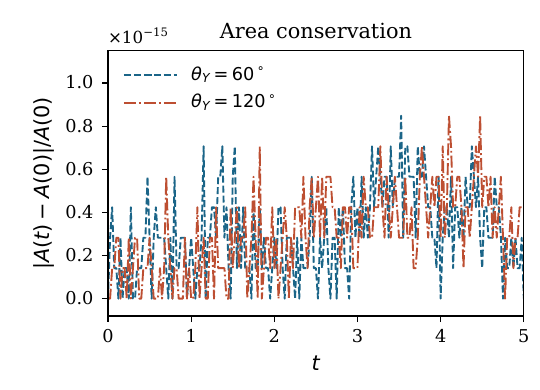}
\caption{}
\end{subfigure}

\begin{subfigure}[b]{0.48\linewidth}
\centering
\includegraphics[width=\linewidth]{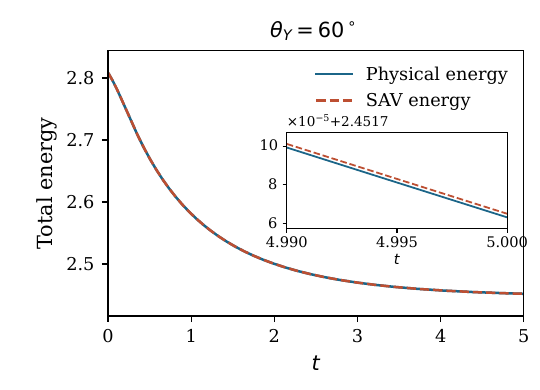}
\caption{}
\end{subfigure}
\begin{subfigure}[b]{0.48\linewidth}
\centering
\includegraphics[width=\linewidth]{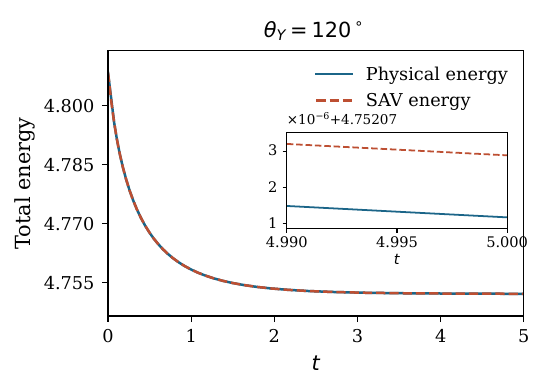}
\caption{}
\end{subfigure}
\vspace{-5pt}
\caption{Evolution of the contact angle (a), the volume error
(b), and the total energy and SAV energy (c)--(d) with $g=1$ for the two-dimensional numerical solutions.}
\label{fig:ghl-gravity-histories}
\end{figure}

Figure~\ref{fig:ghl-gravity-meshes}  shows the computed droplet profiles,
meshes, and velocities for  $g = 0, 1, 2, 10$. The wetting and dewetting case correspond to the static contact angles $\theta_Y = 60^\circ$ and $\theta_Y = 120^\circ$, respectively.  The droplet becomes flatter
as $g$ increases. In both cases, the method maintains good mesh
quality throughout the evolution, including near the moving contact
points, even at $g=10$.

Figure~\ref{fig:ghl-gravity-histories} shows the time evolution of the
contact angle, the scaled volume error for $g=1$, and the discrete
total and SAV energies. The contact angle approaches $\theta_Y$,
and the relative volume error remains on the order of $10^{-16}$,
indicating that the scheme preserves the droplet volume to machine
precision. The physical and SAV energies nearly coincide on the full
scale, and the insets show their late-time differences. In the presence of
gravity, the total energy decreases initially and then reaches a nearly constant
value.

\FloatBarrier

\subsection{Three-dimensional ellipsoid droplet}\label{sec:numerics-3d}
We next consider a three-dimensional droplet problem. The initial
droplet occupies the upper half of an ellipsoid with semiaxes
$(1.5,1,1)$ along the $x_1$, $x_2$, and $y$ directions
\cite{GHL}. We set the initial velocity to zero and use $P_2/P_1$
elements with quadratic geometry. 

\begin{figure}[H]
\centering
\includegraphics[width=0.45\linewidth]{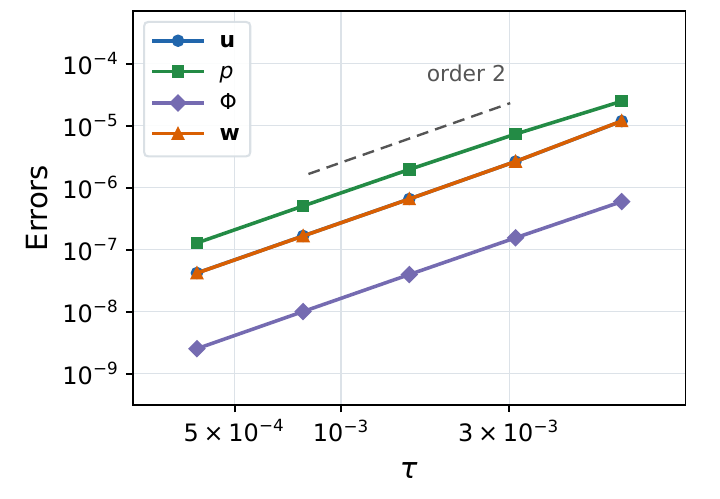}
\includegraphics[width=0.45\linewidth]{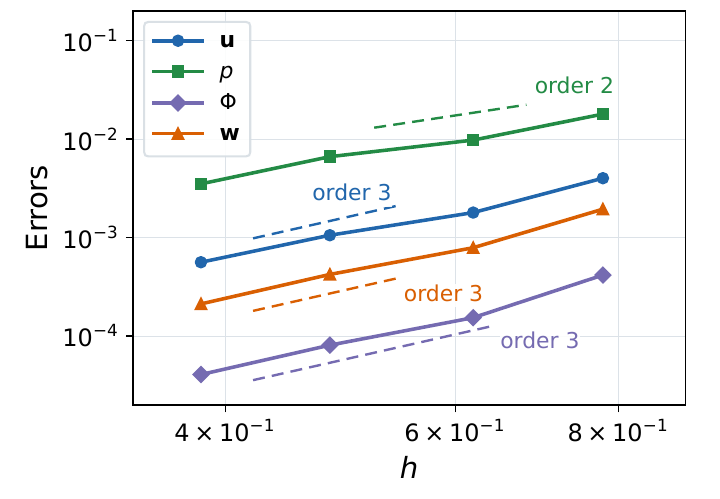}
\vspace{-5pt}
\caption{Temporal discretization errors  (left) and spatial  discretization errors (right) of the three-dimensional numerical solutions.}
\label{fig:time-3d}
\end{figure}

For the convergence tests,
the viscosity coefficient is $\mu=0.1$, the equilibrium contact
angle is $\theta_Y=\pi/2$, $g=\beta_s=\beta_c=0$.
Figure~\ref{fig:time-3d} shows the spatial and temporal convergence
results at $T=0.1$. The left panel
indicates second-order temporal self-convergence for all four
quantities. The right panel shows
spatial convergence rates close to three for $\bm u,\Phi,\bm w$
and two for the pressure. 

\begin{figure}[H]
\centering
\includegraphics[width=.99\linewidth]{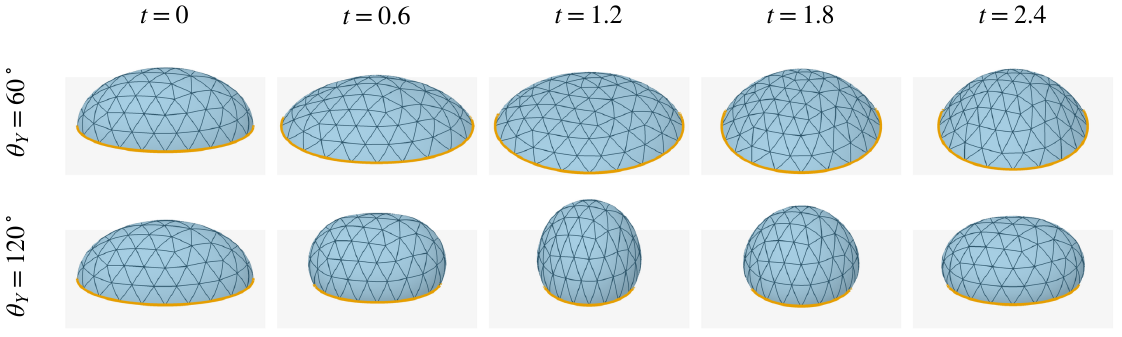}
\caption{Evolution of the ellipsoid droplet profile over one oscillation cycle.}
\label{fig:droplet-3d-shapes}
\end{figure}

We next compare wetting and dewetting from the same initial
droplet under gravitational strength $g=1$. The viscosity,
wall-friction, and contact-line friction coefficients are
$\mu=0.05$, $\beta_s=0.1$, and $\beta_c=0.01$, respectively.
The wetting and dewetting cases correspond to the equilibrium
contact angles $\theta_Y=60^\circ$ and $\theta_Y=120^\circ$,
respectively.
Figure~\ref{fig:droplet-3d-shapes} shows the computed droplet
profiles using a time step size of $\tau=0.005$ up to the final
time $T=2.4$. The droplet spreads and becomes flatter in the
wetting case, whereas it retracts and rises in the dewetting case.
In both cases, the snapshots show good mesh quality near the
moving contact line.

\begin{figure}[H]
\centering
\begin{subfigure}[b]{0.48\linewidth}
\centering
\includegraphics[width=\linewidth]{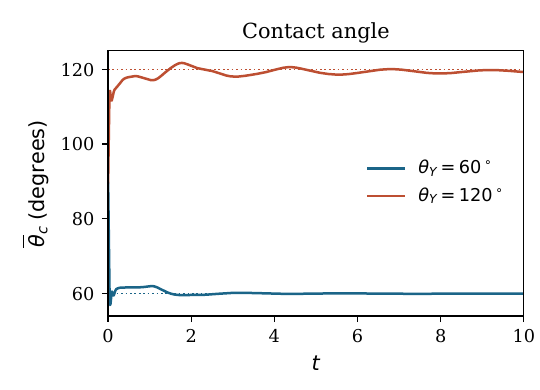}
\caption{}
\end{subfigure}
\begin{subfigure}[b]{0.48\linewidth}
\centering
\includegraphics[width=\linewidth]{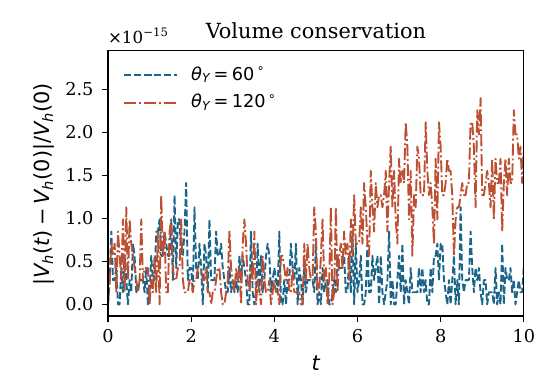}
\caption{}
\end{subfigure}

\begin{subfigure}[b]{0.48\linewidth}
\centering
\includegraphics[width=\linewidth]{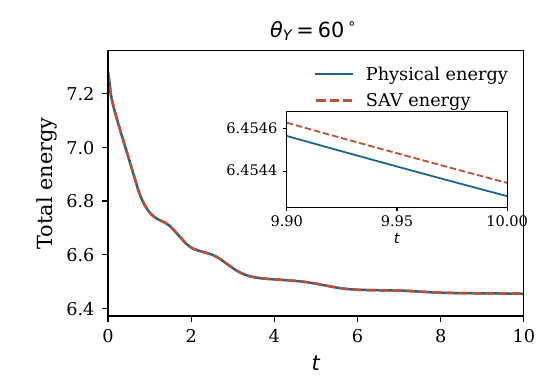}
\caption{}
\end{subfigure}
\begin{subfigure}[b]{0.48\linewidth}
\centering
\includegraphics[width=\linewidth]{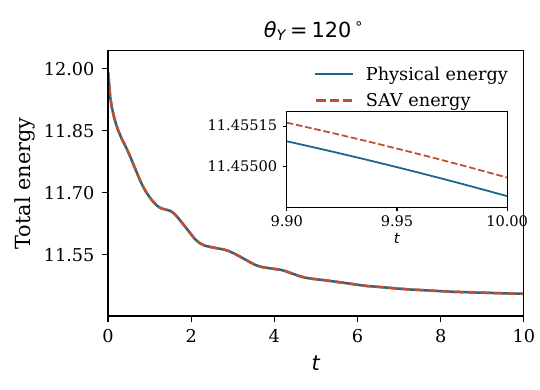}
\caption{}
\end{subfigure}
\caption{Evolution of the contact angle (a), the volume error
(b), and the total energy and SAV energy (c)--(d) with $g=1$ for the three-dimensional numerical solutions.}
\label{fig:droplet-3d-history}
\end{figure}

Figure~\ref{fig:droplet-3d-history} shows the time evolution of
the mean contact angle, the relative volume error, and the
physical and SAV energies. These histories are computed on the
same mesh with a smaller time step size of $\tau=0.002$ up to
$T=10$. The mean contact angles approach their respective
equilibrium values $\theta_Y$, while the relative volume error
stays below $2.7\times10^{-15}$, indicating volume conservation
to machine precision. The SAV energy decreases in both cases.
The physical and SAV energy curves nearly coincide on the scale
of the main panels, and the insets resolve their small late-time
separation.

\section{Conclusions}
We have developed a high-order, structure-preserving
ALE--SAV finite element method for incompressible Navier--Stokes flow
with gravity and moving contact lines. High-order isoparametric spatial
approximation and second-order time stepping are combined with linear
solves at each step. The exactly $L^2$-isometric velocity transfer and
the volume-preserving boundary update obtained through linear solves
make the moving-mesh energy
balance compatible with exact area or volume conservation. We have established
an exact discrete SAV energy-dissipation identity and exact area or
volume conservation on the final mesh. Numerical experiments in two
and three dimensions demonstrate high-order spatial convergence and
second-order temporal convergence, and verify SAV energy dissipation
and area or volume conservation.

\subsection*{Funding}
This work was supported by the National Natural Science Foundation of China (Young Scientists Fund, Category C, Grant No. 1260010446) and the Guangdong Basic and Applied Basic Research Foundation (Grant No. 2026A1515012472).

\FloatBarrier
\clearpage
\bibliographystyle{siamplain-nodoi}
\bibliography{references}
\end{document}